\documentclass[12pt,3p,times]{elsarticle}

\usepackage{amsmath,amsfonts,amssymb}
\usepackage{amsthm}
\usepackage{bm}
\usepackage{graphicx,epstopdf}
\usepackage[position=top,labelfont=normalfont,textfont=normalfont,singlelinecheck=off,justification=raggedright]{subfig}
\usepackage{floatrow}
\usepackage[dvipsnames]{xcolor}
\usepackage{setspace}
\usepackage{enumerate,etaremune}
\usepackage{booktabs}
\usepackage{tabularx,multirow}
\usepackage{framed}
\usepackage{hhline}
\usepackage{url}
\usepackage[unicode,bookmarks=false]{hyperref}
\hypersetup{
  colorlinks,
  citecolor=Blue,linkcolor=Blue,urlcolor=Blue
}
\usepackage[toc,page]{appendix}

\usepackage{lineno}
\usepackage[textsize=footnotesize,linecolor=red,backgroundcolor=red!25,bordercolor=red]
  {todonotes}

\usepackage{tikz}
\usetikzlibrary{shapes.geometric}
\usetikzlibrary{shapes,arrows}
\usetikzlibrary{plotmarks}

\usetikzlibrary{calc}

\usepackage{algorithm}
\usepackage{algorithmicx,algpseudocode}
\usepackage{cases}
\newcommand{\eg}{{\it e.g.}}
\newcommand{\ie}{{\it i.e.}}

\newcommand{\etal}{{\it et al.}}
\newcommand{\tensor}[1]{\bm{#1}}
\newcommand{\stress}{\sigma}
\newcommand{\strain}{\varepsilon}
\newcommand{\tstress}{\tensor{\stress}}
\newcommand{\tstrain}{\tensor{\strain}}

\newcommand{\jump}[1]{\tilde{#1}}

\newcommand{\dd}{\mathrm{d}}
\newcommand{\pd}{\partial}

\newcommand{\el}{\mathrm{e}}

\newcommand{\rn}[1]{\uppercase\expandafter{\romannumeral #1\relax}}
\newcommand{\cn}{N}
\newcommand{\ct}{T}

\DeclareMathOperator{\diver}{\nabla\cdot}
\DeclareMathOperator{\symgrad}{\nabla^{s}}

\newsavebox{\dotbox}

\theoremstyle{remark}
\newtheorem{remark}{Remark}

\definecolor{kaist-blue}{RGB}{20,135,200}
\definecolor{kaist-dark-blue}{RGB}{0,65,145}
\definecolor{kaist-medium-blue}{RGB}{0,65,135}
\definecolor{kaist-light-blue}{RGB}{95,190,235}
\definecolor{kaist-dark-gray}{RGB}{124,124,124}

\newcommand{\revised}[1]{{\color{black} #1}}

\newcolumntype{L}[1]{>{\raggedright\let\newline\\arraybackslash\hspace{0pt}}m{#1}}
\newcolumntype{C}[1]{>{\centering\let\newline\\arraybackslash\hspace{0pt}}m{#1}}
\newcolumntype{R}[1]{>{\raggedleft\let\newline\\arraybackslash\hspace{0pt}}m{#1}}

\allowdisplaybreaks

\biboptions{sort&compress,square,comma,numbers}
\AtBeginDocument{\hypersetup{citecolor=MidnightBlue,linkcolor=MidnightBlue,urlcolor=MidnightBlue}}

\usepackage{etoolbox}
\makeatletter
\patchcmd{\ps@pprintTitle}% <cmd>
{Preprint submitted to}% <search>
{~}% <replace>
{}{}% <succes><failure>
\makeatother
\begin{document}

\begin{frontmatter}

\title{A barrier method for frictional contact on embedded interfaces}

\author[KAIST]{Yidong Zhao\fnref{equal}}
\author[KAIST]{Jinhyun Choo\fnref{equal}\corref{corr}}
\ead{jinhyun.choo@kaist.ac.kr}
\author[HKU]{Yupeng Jiang}
\author[UCLA]{Minchen Li}
\author[UCLA]{Chenfanfu Jiang}
\author[UCB]{Kenichi Soga}

\cortext[corr]{Corresponding Author}
\fntext[equal]{These authors contributed equally to this work}

\address[KAIST]{Department of Civil and Environmental Engineering, KAIST, South Korea}
\address[HKU]{Department of Civil Engineering, The University of Hong Kong, Hong Kong}
\address[UCLA]{Department of Mathematics, University of California, Los Angeles, United States}
\address[UCB]{Department of Civil and Environmental Engineering, University of California, Berkeley, United States}

\journal{~}

\begin{abstract}
We present a barrier method for treating frictional contact on interfaces embedded in finite elements.
The barrier treatment has several attractive features, including:
(i) it does not introduce any additional degrees of freedom or iterative steps,
(ii) it is free of inter-penetration,
(iii) it avoids an ill-conditioned matrix system,
and (iv) it allows one to control the solution accuracy directly.
We derive the contact pressure from a smooth barrier energy function that is designed to satisfy the non-penetration constraint.
Likewise, we make use of a smoothed friction law in which the stick--slip transition is described by a continuous function of the slip displacement.
We discretize the formulation using the extended finite element method to embed interfaces inside elements, and devise an averaged surface integration scheme that effectively provides stable solutions without traction oscillations.
Subsequently, we develop a way to tailor the parameters of the barrier method to embedded interfaces, such that the method can be used without parameter tuning.
We verify and investigate the proposed method through numerical examples with varied levels of complexity.
The numerical results demonstrate that the proposed method is remarkably robust for challenging frictional contact problems, while requiring low cost comparable to that of the penalty method.
\end{abstract}

\begin{keyword}
Barrier method \sep
Contact algorithm\sep
Frictional sliding \sep
Embedded interfaces \sep
Extended finite element method
\end{keyword}

\end{frontmatter}

% \linenumbers

% SECTION
% ==============================================================================
\section{Introduction}
Frictional interfaces inside solids appear in a wide range of problems in science and engineering.
Examples in civil and mechanical engineering applications include cracks, slip surfaces, and soil--structure interfaces.
\revised{Mathematically, frictional contact is a constrained optimization problem, in which the constraints emanate from non-penetration of two contacting surfaces and frictional resistance in the surfaces.
Therefore, a numerical strategy for the contact problem should deal with two related aspects: (i) the algorithm to impose the constraints, and (ii) the discretization of the constraint variables.} 

In modern finite element analysis, interfaces are often embedded in elements (\ie~allowed to pass through the interior of elements) to simplify the meshing procedure dramatically.
When the interfaces are modeled as sharp (lower-dimensional) discontinuities, they are embedded by enriching basis functions either locally or globally.
Representative examples of locally or globally enriched finite element methods are the assumed enhanced strain (AES) method~\cite{simo1990class} and the extended finite element method (XFEM)~\cite{moes1999finite}, respectively.
At the expense of the simplified mesh, however, these methods pose a new challenge of treating frictional contact on embedded interfaces.

Over the last couple of decades, a variety of methods have been studied for modeling frictional interfaces embedded in finite elements.
Dolbow \etal~\cite{dolbow2001extended} were the first to propose a method for enforcing frictional contact on interfaces embedded with XFEM. Their method employs an iterative scheme called LATIN to solve a variational equation formulated in terms of the total displacements of the interface surfaces.
Yet the convergence behavior and accuracy of the LATIN iterative scheme appeared to have much room for improvement.
As such, a number of subsequent studies have developed and investigated various types of methods for a more efficient and robust treatment of embedded frictional interfaces (\eg~\cite{khoei2007enriched,liu2008contact,nistor2009xfem}).
A notable example is the work of Liu and Borja~\cite{liu2008contact}, where the variational equation is reformulated based on the relative displacements of the discontinuity surfaces and the contact constraints are treated by the classical penalty method.
Other classical methods for handling frictional contact, such as the Lagrange multiplier method and the augmented Lagrangian method, have also been applied and tailored to embedded interfaces (\eg~\cite{elguedj2007mixed,bechet2009stable,liu2010stabilized}).
\revised{At the same time, some researchers have pursued alternative discretization strategies for the contact variables, instead of discretizing them in a way similar to the classical node-to-surface method in computational contact mechanics~\cite{laursen2003computational,wriggers2006computational}.
A remarkable example is the mortared finite element formulation proposed by Kim \etal~\cite{kim2007mortared}, where the contact variables are discretized on an intermediate mortar surface independent from the bulk mesh.
In that formulation, the displacement fields in the bulk material and the mortar surface are weakly coupled via a Lagrange multiplier, and the contact constraints are handled by the penalty method.}

Nevertheless, it remains an unresolved challenge as to how to handle frictional contact on embedded interfaces in a both computationally robust and efficient manner.
The Lagrange multiplier method provides exact enforcement of the contact constraints.
However, it introduces an additional field variable which not only increases the number of degrees of freedom but also gives rise to a stability issue in mixed discretization~\cite{ji2004strategies,bechet2009stable,liu2010stabilized}.
As a practical means, the penalty method is often used because it retains the number of degrees of freedom and can be used easily.
However, this method unavoidably permits inter-penetration (over-closure) to some extent and the degree of inter-penetration can only be reduced indirectly by enlarging the penalty stiffness parameter.
Unfortunately, too large a penalty parameter makes the resultant matrix system highly ill-conditioned.
So the determination of an appropriate penalty parameter is often a non-trivial issue.
\revised{The augmented Lagrangian method combines ideas from the Lagrange multiplier and the penalty methods to alleviate the drawbacks of each method.
However, it is still subjected to the stability problem (when the multiplier is discretized independently as in Alart and Curnier~\cite{alart1991mixed}) or requires additional iterative steps to calculate the multiplier (when the multiplier is not discretized separately as in Simo and Laursen~\cite{simo1992augmented}).}
Recent studies have thus developed alternative and more advanced methods such as the weighted Nitsche method~\cite{annavarapu2013nitsche,annavarapu2014nitsche,annavarapu2015weighted} or explored diffuse approaches for algorithm-free modeling of frictional interfaces~\cite{fei2020phase-a,fei2020phase-b,fei2021double,fei2021phase}.
Yet none of the existing methods is considered optimally robust and efficient for modeling embedded interfaces with frictional contact.

The barrier method -- a simple yet robust algorithm for constrained optimization -- has also been used for imposing the contact constraints on external domain boundaries in multi-body dynamics~\cite{kane1999finite,pandolfi2002time,li2020incremental,Li2021CIPC,Ferguson:2021:RigidIPC,Lan2021MIPC}.
In essence, the method replaces the constraint of a constrained optimization problem by augmenting a barrier function to the objective functional, so that one can solve the problem using an algorithm for unconstrained optimization such as Newton's method.
In the context of contact mechanics, the non-penetration constraint can be approximated as a barrier energy function augmented to the potential energy functional.
The particular barrier energy function used in Li \etal~\cite{li2020incremental} is a product of the following two: (i) a $C^{2}$-continuous barrier function of the unsigned distance between domain boundaries,
and (ii) a scalar scaling parameter that converts the value of the barrier function to an energy.
Importantly, the barrier function has a free parameter that controls the size of the barrier region, whereas the scaling parameter can be tuned to avoid ill-conditioning of the resultant matrix system.
So this barrier method has several advantages over the classical methods, including:
(i) it does not introduce any additional degrees of freedom or iterative steps,
(ii) it is free of inter-penetration,
(iii) it avoids an ill-conditioned matrix system,
and (iv) it allows one to control the solution accuracy directly.
\revised{It is noted that although the barrier method results in a very small gap between two contacting surfaces, such a gap indeed exists in real world due to asperities at microscopic scales.
In this regard, a gap in the barrier method is more physically realistic than an inter-penetration in the penalty method.
Further, whenever the crack opening displacement is used to calculate physical properties such as the hydraulic conductivity (\eg~\cite{liu2017stabilized,choo2018cracking,liu2020modeling}), a gap must be preferred to an inter-penetration to ensure non-negative values of the properties.}

While the features of the barrier method are also attractive for embedded interfaces, the existing barrier method is incompatible with interfaces embedded in finite elements.
A few critical reasons are as follows.
First, the majority of embedded interfaces are closed initially, but the barrier method does not allow one to initialize perfectly closed surfaces because the barrier energy becomes infinity as the distance between the two surfaces approaches zero.
An intuitive way to address this issue is to introduce a new parameter representing the initial ``gap'' of a closed interface; however, it is unclear how to determine the value of this new parameter.
Second, it is uncertain as to how to determine the scalar scaling parameter for an embedded interface under quasi-static conditions, because the existing way to determine the parameter assumes dynamic contact between multiple bodies.
Third, in embedded finite element methods, numerical solutions may show severe oscillations in traction fields when the enriched degrees of freedom pose strong constraints.

In this work, we develop the first barrier method for embedded interfaces with frictional contact, addressing the aforementioned challenges in a physically meaningful and numerically efficient manner.
The work proceeds as follows.
After formulating a variational equation for a solid with internal discontinuities, we employ a barrier energy function to derive the contact pressure such that the non-penetration constraint is always satisfied.
Likewise, we utilize a smoothed friction law in which the stick--slip transition is a continuous function of the slip displacement.
We then discretize the formulation using XFEM to embed interfaces inside finite elements.
In doing so, we devise a surface integration scheme for XFEM that can simply suppress oscillations in traction fields.
Subsequently, we develop a way to tailor the parameters of the barrier method to embedded interface, such that the resulting method can be used without parameter tuning.
Through numerical examples with varied levels of complexity, we verify and investigate the proposed method for handling embedded frictional interfaces, and discuss its performance compared with that of the penalty method which has similar computational cost. To limit the scope of the present work, we shall focus on single and stationary interfaces \revised{in two dimensions} throughout.
\revised{It is noted that the present scope is comparable to that of previous work on frictional contact on embedded interfaces (\eg~\cite{annavarapu2014nitsche}).}

% SECTION
% ==============================================================================
\section{Formulation}

In this section, we first formulate a variational equation for a solid with internal discontinuities, which are subjected to constraints on the contact condition and frictional sliding.
We then handle the constraints employing a smooth barrier energy function and a smoothed friction law.
Without loss of generality, we assume the quasi-static condition, absence of body force, and infinitesimal deformation.

\subsection{Problem statement}
Consider a domain $\Omega\in\mathbb{R}^{\rm dim}$ delimited by boundary $\pd\Omega$ with outward normal vector $\bm{\upsilon}$.
The boundary is suitably decomposed into the displacement (Dirichlet) boundary $\pd_{u}\Omega$ where the displacement vector is prescribed as $\bar{\bm{u}}$, and the traction (Neumann) boundary $\pd_{t}\Omega$ where the traction vector is prescribed as $\bar{\bm{t}}$, such that $\pd\Omega=\overline{\pd_{u}\Omega\cup\pd_{t}\Omega}$ and $\emptyset=\pd_{u}\Omega\cap\pd_{t}\Omega$.
The domain contains a set of embedded, lower-dimensional interfaces $\Gamma$, across which the displacement field $\bm{u}(\bm{x})$ can be discontinuous ($\bm{x}$ denotes the position vector).
Each interface segment has two surfaces, and the positive and negative surfaces are denoted by $\Gamma_{+}$ and $\Gamma_{-}$, respectively.
In that follows, we shall use the subscripts $(\cdot)_{+}$ and $(\cdot)_{-}$ to denote quantities pertaining to the positive and negative sides.
For example, the two sides of the domain are denoted by $\Omega_{+}$ and $\Omega_{-}$.
We denote by $\bm{n}\equiv\bm{n}_{+}=-\bm{n}_{-}$ the unit normal vector to the interface pointing toward $\Gamma_{+}$ and by $\bm{m}$ the unit tangent vector pointing toward the slip direction.

Let $\jump{\bm{u}}:=\bm{u}_{+}-\bm{u}_{-}$ denote the relative displacement (displacement jump) across the interface surfaces.
The gap distance between the two surfaces is given by
\begin{equation}
  u_{\cn} = \jump{\bm{u}}\cdot\bm{n}.
\end{equation}
The non-penetration constraint then can be written as
\begin{equation}
  u_{\cn} \geq 0.
  \label{eq:non-penetration-constraint}
\end{equation}

Let $\bm{t}\equiv\bm{t}_{+}=-\bm{t}_{-}$ denote the traction vector acting on the contact surfaces. 
The traction vector is decomposed into its normal and tangential components as
\begin{equation}
  \bm{t} = -p_{\cn}\bm{n} + \bm{t}_{\ct},
  \label{eq:contact-traction-decomposition}
\end{equation}
where $p_{\cn}$ is the contact pressure, and $\bm{t}_{\ct}$ is the tangential traction vector.
It is noted that the contact pressure should be non-negative, \ie
\begin{equation}
  p_{\cn} \geq 0.
  \label{eq:contact-pressure-constraint}
\end{equation}

The frictional sliding behavior is constrained by a friction law.
In this work, we consider the standard Coulomb friction law, which can be written as
\begin{linenomath}
\begin{align}
  \tau \leq \mu p_{\cn},
  \label{eq:friction-law}
\end{align}
\end{linenomath}
where $\tau:=\|\bm{t}_{T}\|$ is the tangential stress (resolved shear stress), and $\mu$ is the friction coefficient.
Conventionally, the constraint imposed by the friction law has been cast into a yield function.
However, we shall not follow the conventional method and introduce a new approach later in this section.

For simplicity and without loss of generality, we assume that the bulk material is linear elastic.
The constitutive behavior of the bulk material can then be written as
\begin{equation}
  \tstress = \mathbb{C}^{\el}:\tstrain,
\end{equation}
where $\tstress$ is the Cauchy stress tensor, $\tstrain:=\symgrad\bm{u}$ is the infinitesimal strain tensor defined as the symmetric gradient of the displacement field, and $\mathbb{C}^{\el}$ is the fourth-order elasticity stiffness tensor.

The strong form of the boundary-value problem can be stated as follows: find the displacement field $\bm{u}$ such that
\begin{linenomath}
\begin{align}
    \diver\bm{\sigma} = \bm{0} \;\; &\text{in}\;\; \Omega\setminus\Gamma,
    \label{eq:linear-momentum} \\
    \bm{u} = \bar{\bm{u}} \;\; &\text{on}\;\; \pd_{u}\Omega,
    \label{eq:displacement-bc} \\
    \bm{\upsilon}\cdot\tstress = \bar{\bm{t}} \;\; &\text{on}\;\; \pd_{t}\Omega,
    \label{eq:traction-bc} \\
    \bm{n}\cdot\tstress = \bm{t}_{-} \;\; &\text{on}\;\; \Gamma_{-}, \\
    -\bm{n}\cdot\tstress = \bm{t}_{+} \;\; &\text{on}\;\; \Gamma_{+}.
\end{align}
\end{linenomath}

To develop the weak form of the problem, we define the spaces of trial functions and variations, respectively, as
\begin{linenomath}
\begin{align}
    \bm{\mathcal{U}} &:= \{\bm{u}\,|\,\bm{u}\in H^{1}(\Omega),\; \bm{u}=\bar{\bm{u}}\; \text{on}\;  \pd_{u}\Omega\}, \\
    \bm{\mathcal{V}} &:= \{\bm{\eta}\,|\,\bm{\eta}\in H^{1}(\Omega),\; \bm{\eta}=\bm{0}\; \text{on}\;  \pd_{u}\Omega\},
\end{align}
\end{linenomath}
where $H^{1}$ is the Sobolev space of order one.
Through the standard procedure, we arrive at the following variational equation
\begin{linenomath}
\begin{align}
  \int_{\Omega\setminus\Gamma} \symgrad{\bm{\eta}}:\tstress\, \dd V
  + \int_{\Gamma} \jump{\bm{\eta}}\cdot\bm{t}\, \dd A
  - \int_{\pd_{t}\Omega} \bm{\eta}\cdot\bar{\bm{t}}\, \dd A
  = 0,
  \label{eq:variational-form}
\end{align}
\end{linenomath}
where $\jump{\bm{\eta}}:=\bm{\eta}_{+} - \bm{\eta}_{-}$.

The focus of this work is on how to efficiently handle the second term in Eq.~\eqref{eq:variational-form} -- the virtual work done by the contact traction vector -- subjected to the contact constraints.
For this purpose, we decompose the contact surface integral into its normal and tangential components by substituting Eq.~\eqref{eq:contact-traction-decomposition}, \ie
\begin{equation}
    \int_{\Gamma} \jump{\bm{\eta}}\cdot\bm{t}\, \dd A
    = \int_{\Gamma} \jump{\bm{\eta}}\cdot(-p_{\cn}\bm{n})\, \dd A
    + \int_{\Gamma} \jump{\bm{\eta}}\cdot\bm{t}_{\ct}\, \dd A.
    \label{eq:contact-integral}
\end{equation}
In what follows, we first derive $p_{\cn}$ from a barrier energy function to address the contact normal behavior, and then use a smoothed friction law to evaluate $\bm{t}_{\ct}$ in the sliding behavior.

\subsection{Barrier method}
In this work, we make use of a barrier method to treat the contact normal behavior subjected to the non-penetration constraint.
While one can apply the barrier method in a purely mathematical manner as in Li \etal~\cite{li2020incremental}, here we do it from a physical point of view to provide another perspective of the method.
The physical perspective would be useful to better understand the barrier method, like how the spring interpretation of the penalty method does. 

\begin{figure}[h!]
    \centering
    \includegraphics[width=0.5\textwidth]{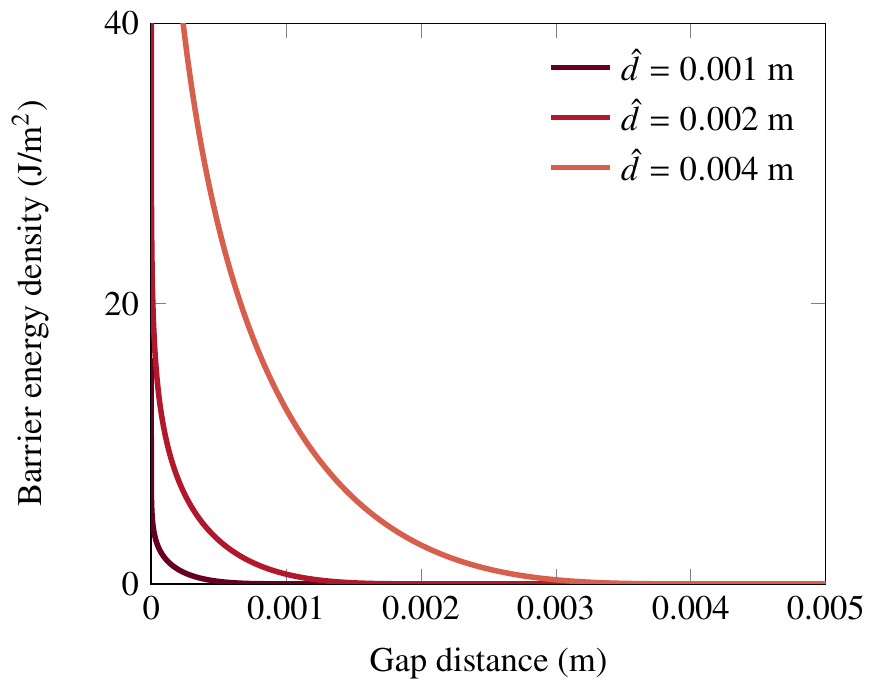}
    \caption{Variation of the barrier energy density function, $B(u_{\cn})$, with $\hat{d}$. ($\kappa = 1$ MPa/m)}
    \label{fig:barrier-functions}
\end{figure}

To begin, let us introduce a elastic barrier between the two faces of a discontinuity.
The elastic energy stored in the barrier can be represented by a barrier energy density function taking the gap distance as an argument.
In this work, we adapt the smooth barrier energy function proposed by Li \etal~\cite{li2020incremental} to embedded interfaces, which gives
\begin{linenomath}
\begin{align}
  B(u_{\cn}) := 
  \begin{cases}
  -\kappa(u_{\cn} - \hat{d})^{2}\ln\left(\dfrac{u_{\cn}}{\hat{d}}\right) & \text{if}\;\; 0<u_{\cn}<\hat{d}, \\
  0 & \text{if}\;\; u_{\cn}\geq\hat{d}.
  \end{cases}
  \label{eq:barrier-function}
\end{align}
\end{linenomath}
Here, $\hat{d}$ is the value of $u_{\cn}$ above which the value of the barrier function is zero.
In other words, the value of $\hat{d}$ defines the maximum value of gap distance at which the two surfaces are considered ``in contact.''
In this sense, $\hat{d}$ can be interpreted as the thickness of the barrier.
Also, $\kappa>0$ is a scalar parameter having the unit of pressure per length, which is introduced to let the unit of $B(u_{\cn})$ be energy per area.
Later in this section, we will show that $\kappa$ controls the stiffness of the barrier for a given $\hat{d}$. 
Figure~\ref{fig:barrier-functions} illustrates how the values of the barrier energy density function vary with $\hat{d}$.

\begin{figure}[h!]
    \centering
    \includegraphics[width=0.5\textwidth]{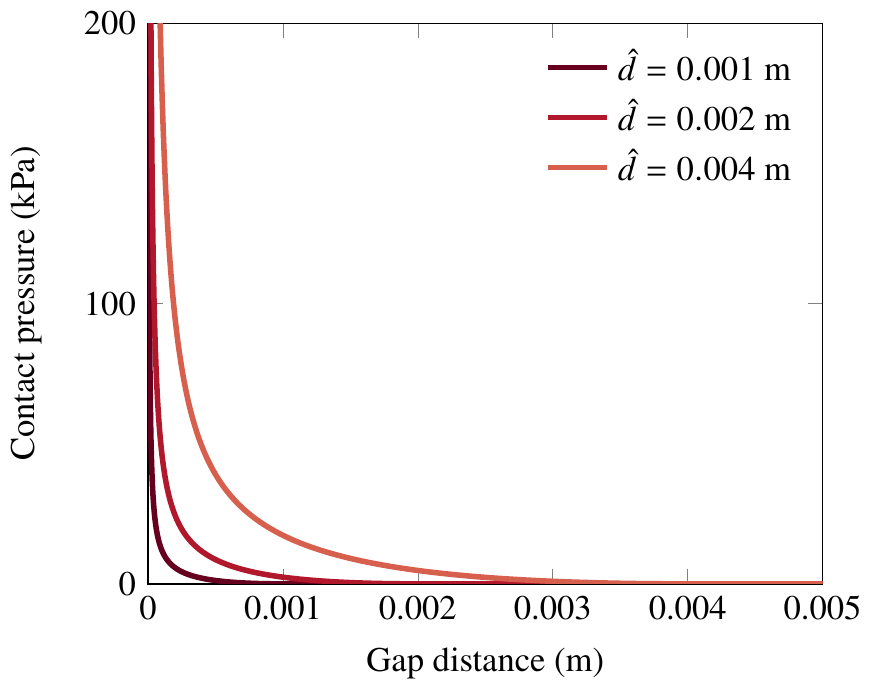}
    \caption{Variations of the contact pressure, $p_{\cn}$, with the gap distance, $u_{\cn}$, derived from the barrier energy functions in Fig.~\ref{fig:barrier-functions}.}
    \label{fig:contact-pressure}
\end{figure}

The contact pressure can then be derived based on energy conjugacy, \ie
\begin{linenomath}
\begin{align}
  p_{N} := 
  -\dfrac{\pd B(u_{\cn})}{\pd u_{\cn}} =
  \begin{cases}
  \kappa(u_{\cn} - \hat{d})\left[2\ln\left(\dfrac{u_{\cn}}{\hat{d}}\right) - \dfrac{\hat{d}}{u_{\cn}} + 1\right] & \text{if}\;\; 0<u_{\cn}<\hat{d}, \\
  0 & \text{if}\;\; u_{\cn} \geq \hat{d}.
  \end{cases}
  \label{eq:contact-pressure}
\end{align}
\end{linenomath}
Figure~\ref{fig:contact-pressure} shows the relations between the contact pressure and the gap distance derived from the barrier functions in Fig.~\ref{fig:barrier-functions}.
One can see that the contact pressure becomes positive when the gap distance becomes smaller than $\hat{d}$ and it increases toward infinity as the gap approaches zero.
The contact pressure thus prevents the gap to become zero, not to mention negative, ensuring satisfaction of the non-penetration constraint, Eq.~\eqref{eq:non-penetration-constraint}.
It is also noted that the contact pressure is always non-negative, satisfying Eq.~\eqref{eq:contact-pressure-constraint}.

Further, we can gain insight into the physical meaning of $\kappa$ by calculating the tangent of the (negative) contact pressure with respect to the gap distance, which can be regarded as the normal stiffness of the interface.
The tangent is given by
\begin{linenomath}
\begin{align}
    k_{\cn} := 
    \dfrac{\pd (-p_{\cn})}{\pd u_{\cn}} =
    \begin{cases}
    -2\kappa\ln\left(\dfrac{u_{\cn}}{\hat{d}}\right) - \kappa\dfrac{(u_{\cn}-\hat{d})(3u_{\cn} + \hat{d})}{u_{\cn}^{2}} & \text{if}\;\; 0<u_{\cn}<\hat{d}, \\
    0 & \text{if}\;\; u_{\cn} \geq \hat{d}.
    \end{cases}
    \label{eq:barrier-stiffness}
\end{align}
\end{linenomath}
The above equation shows that the value of $\kappa$ directly controls the normal stiffness of the barrier-treated interface.
For this reason, hereafter we shall refer to $\kappa$ as the barrier stiffness parameter.
It is noted that this parameter affects not only the physical behavior of the barrier but also the matrix condition of the discretized problem.
Later in Section~\ref{sec:parameterization}, we will explain how to determine this parameter for an embedded interface such that it results in a physically meaningful and numerically efficient formulation.

\smallskip
\begin{remark}
    The barrier method can be easily implemented and utilized with a trivial modification of any existing penalty-based code.
    Concretely, one only needs to replace the contact pressure--penetration relationship in the penalty method with the contact pressure--gap distance relationship in the barrier method, Eq.~\eqref{eq:contact-pressure}.
    So the computational cost of the barrier method is practically the same as that of the penalty method.
\end{remark}

\smallskip
\begin{remark}
    The form of the barrier function~\eqref{eq:barrier-function}, proposed in Li \etal~\cite{li2020incremental}, has two features that distinguish itself from the standard logarithmic barrier function in optimization~\cite{boyd2004convex}: (i) it is truncated at $\hat{d}$, and (ii) a quadratic term is multiplied to a logarithmic term.
    The first feature allows us not only to avoid the barrier energy to be negative but also to localize the barrier energy such that the method can be efficiently utilized for multiple pairs of contacting surfaces. 
    The second feature makes the barrier energy function $C^{2}$-continuous despite the truncation.
    The $C^{2}$-continuity is critical to solving the problem using a gradient-based algorithm such as Newton's method.
\end{remark}

\smallskip
\begin{remark}
    The use of the barrier method has made the non-penetration constraint more strict (\ie~$u_{\cn}>0$) since $p_{\cn}(u_{\cn})\rightarrow\infty$ as $u_{\cn}\rightarrow 0$.
    While this change may be viewed as a source of numerical inaccuracy, it can also be physically justified in that no real surfaces can be perfectly closed together due to asperities at smaller scales.
    In other words, surfaces can only be closed at the scales of their asperities.
    A similar argument can be found from Wriggers \etal~\cite{wriggers1990finite}, among others.
    \label{rem:impossible-zero-gap}
\end{remark}

\smallskip
\begin{remark}
    One important difference between embedded interface problems and multi-body contact problems is that embedded interfaces are often closed in the initial condition.
    For the reason explained in Remark~\ref{rem:impossible-zero-gap}, however, a perfectly closed discontinuity cannot be initialized in the barrier method.
    To approximate an initially closed discontinuity, we need a new parameter that represents the initial ``gap'' distance of a closed interface.
    We will thoroughly discuss this new parameter later in Section~\ref{sec:parameterization}, because its appropriate value can be determined after discretization.
\end{remark}

\revised{
\smallskip
\begin{remark}
    While the value of $\hat{d}$ is considered constant in the barrier method, it can be updated through algorithmic iterations -- like interior point methods~\cite{boyd2004convex} -- to attain a more accurate solution.
    At the expense of improved accuracy, however, such an update makes the method more complicated. 
    As will be demonstrated later, the barrier method gives sufficiently good solutions for engineering purposes, provided that the value of $\hat{d}$ is small enough.
    Further, because $\hat{d}$ is the upper bound of the tolerable separation gap between contacting surfaces, it would often be more desirable to specify this upper bound explicitly rather than letting it depend implicitly on numerical parameters as in interior point methods.
    Therefore, here we treat the barrier parameter (and hence the barrier function) constant such that the resultant method can be utilized as efficiently and conveniently as the classic penalty method.
\end{remark}
}

\subsection{Smoothed friction law} 
We now shift our attention to the frictional sliding along the interface.
While the sliding process can be described independently from the foregoing barrier treatment,
here we utilize a smoothed friction law that models the stick--slip transition behavior as a continuous function of the slip displacement~\cite{li2020incremental}.
The use of the smoothed friction law is motivated for two reasons.
First, this smoothing approach is consistent with the barrier treatment where the contact pressure varies smoothly with the gap distance.
Second, it provides significant robustness for problems in which stick and slip conditions are mixed, because the tangential traction increases mildly in the beginning of slip, in contrast to a jump-like increase in the classic yield-function approach.

\begin{figure}[h!]
    \centering
    \includegraphics[width=0.5\textwidth]{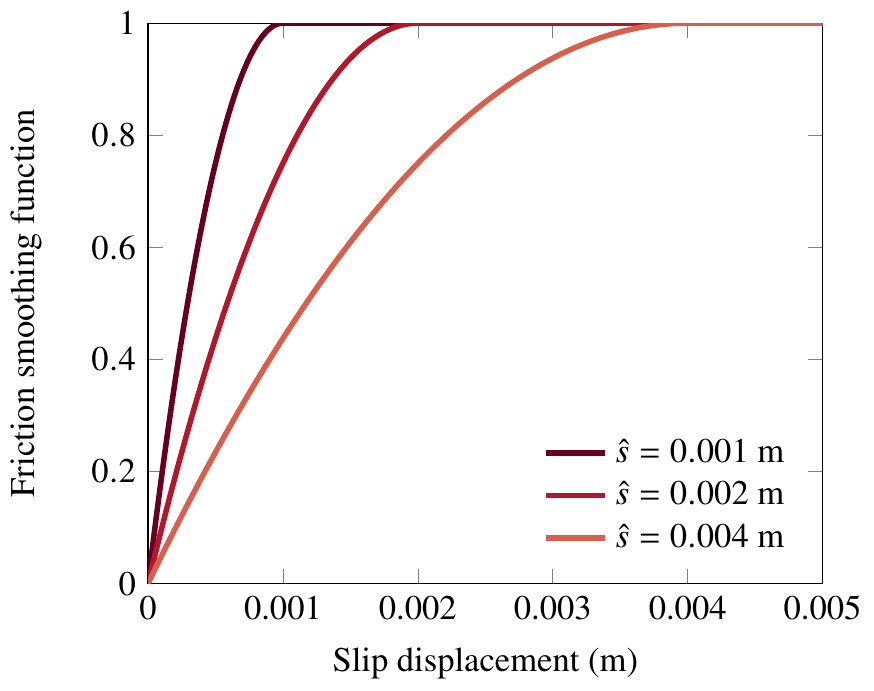}
    \caption{Variation of the friction smoothing function, $m(u_{\ct})$, with $\hat{s}$.}
    \label{fig:smoothed-friction}
\end{figure}

To begin, we introduce a friction smoothing function $m(u_T)\in[0,1]$ which continuously varies with the slip displacement magnitude $u_T := \jump{\bm{u}} \cdot \bm{m}$.
Similar to how we have obtained the barrier energy function~\eqref{eq:barrier-function}, we modify the smoothing function in Li \etal~\cite{li2020incremental} -- originally proposed for multi-body contact dynamics -- for an embedded interface under quasi-static condition.
The modified version reads
\begin{equation}
    m(u_T) := 
    \begin{cases}
    -\dfrac{u_T^2}{\hat{s}^2} + \dfrac{2 |u_T|}{\hat{s}} & \text{if}\; |u_T| < \hat{s}, \\
    1 & \text{if}\; |u_T| \geq \hat{s}.
    \end{cases}
    \label{eq:smoothed-friction-magnitude}
\end{equation}
Here, $\hat{s}$ is a parameter defining the amount of slip displacement at which the kinetic (dynamic) friction is mobilized.
In other words, the function is designed to allow for a small amount of slip in the static friction regime, which is often called microslip.
In this regard, $\hat{s}$ may be called the maximum microslip displacement.
Figure~\ref{fig:smoothed-friction} shows how the friction smoothing function varies with $\hat{s}$.

We then multiply the smoothing function to the tangential stress constrained by the friction law.
When the Coulomb friction law~\eqref{eq:friction-law} is used, the smoothed friction law is
\begin{equation}
    \tau = m(u_T)\mu p_N.
\end{equation}
Accordingly, the tangential traction vector is given by
\begin{equation}
    \bm{t}_T = \tau\bm{m} = m(u_T)\mu p_N\bm{m}.
    \label{eq:smoothed-tangential-traction}
\end{equation}
Notably, Eq.~\eqref{eq:smoothed-tangential-traction} is used to calculate the tangential traction without any yield function.

The main advantage of this smoothed friction law is that the derivative of $\bm{t}_T$ with respect to $u_T$ is a continuous function of $u_T$, because $m(u_T)$ is $C^{1}$-continuous. 
This is in contrast with the classic yield function approach whereby the derivative changes discontinuously between zero and a non-zero value during a stick--slip transition.
Such a continuous change in the derivative of the smoothed friction law provides significant robustness when a stick--slip transition takes place during a gradient-based solution stage such as a Newton iteration.
This aspect will be demonstrated later through a numerical example.

\smallskip
\begin{remark}
    Although we have distinguished the parameter $\hat{s}$ from $\hat{d}$ for its mathematical definition, in practice one can simply assign the same value for these two parameters. 
    According to our experience, any value of $\hat{s}$ works well as long as it is less than the maximum slip displacement.
    Because a reasonable value of $\hat{d}$ must be much smaller than the maximum slip displacement, it is feasible and practical to let $\hat{s}=\hat{d}$.
    By doing so, only the value of $\hat{d}$ needs to be prescribed at one's discretion.
\end{remark}

\smallskip
\begin{remark}
    Like the gap displacement ($u_{\cn} < \hat{d}$) arising in the barrier treatment, the slip displacement in the static friction regime ($u_{\ct} < \hat{s}$) can also be justified from a physical viewpoint. 
    Due to asperities, real contacting surfaces often exhibit a measurable amount of microslip displacement in the static friction regime (\eg~\cite{sandeep2019experimental}).
    Indeed, some researchers have intentionally incorporated microslip displacement using a penalty-like formulation, see, \eg~Wriggers \etal~\cite{wriggers1990finite}.
    In this regard, the smoothed friction law may even be viewed more realistic than the conventional method.
    We also note that the smoothed friction law has two key advantages over the penalty-like model accounting for microslip: (i) it allows us to directly control the maximum microslip displacement via $\hat{s}$, and (ii) it is numerically more robust thanks to its smooth variation with the slip displacement.
    \label{rem:microslip}
\end{remark}

% SECTION
% ==============================================================================
\section{Discretization and algorithms}

In this section, we discretize the formulation using finite elements with embedded discontinuities.
Broadly speaking, there are two classes of finite element methods with embedded discontinuities: (i) the AES method whereby the degrees of freedom of cracked elements are enriched locally~\cite{simo1990class}, (ii) the XFEM approach whereby the degrees of freedom are enriched globally~\cite{moes1999finite}.
Both classes of methods rely commonly on enrichment of finite elements, while the specific kinematic modes of enriched basis functions differ by the particular type of method.
Given that the XFEM enrichment introduces complete kinematic modes for cracked elements (see, \eg~Cusini \etal~\cite{cusini2021simulation} for a relevant discussion), here we choose XFEM as our particular means to study the performance of the barrier method.
We note that the barrier formulation can also be well discretized by any other types of embedded finite element methods, and we believe that the conclusions drawn from the present work may be applied equally well to interfaces embedded by other types of methods.

\subsection{Enrichment of finite elements}
Enrichment of finite elements builds on the idea of decomposing the solution field into its continuous part -- which can be approximated by the standard continuous finite elements -- and its discontinuous part -- which can be approximated by enriching the degrees of freedom of cracked elements.
For the displacement field $\bm{u}(\bm{x})$ passing through the interfaces 
$\Gamma$, let $\bar{\bm{u}}(\bm{x})$ and $\jump{\bm{u}}(\bm{x})$ denote its continuous and discontinuous parts, respectively.
Then the decomposition of the displacement field can be written as~\cite{borja2008assumed}
\begin{equation}
    \bm{u} (\bm{x}) = \bar{\bm{u}}(\bm{x}) + M_{\Gamma}(\bm{x}) \jump{\bm{u}} (\bm{x}).
    \label{eq:displacement-decomposition}
\end{equation}
Here, $M_{\Gamma}(\bm{x})$ is a scalar function for enrichment, given by
\begin{equation}
    M_{\Gamma}(\bm{x}) = H_{\Gamma} (\bm{x}) - f^h (\bm{x}),
\end{equation}
where $H_{\Gamma} (\bm{x})$ is the Heaviside function, which is discontinuous on the surface $\Gamma$ as
\begin{equation}
    H_{\Gamma} (\bm{x}) = 
    \begin{cases}
        1 & \text{if}\;\; \bm{x} \in \Omega_{+} \\
        0 & \text{if}\;\; \bm{x} \in \Omega_{-}\,,
    \end{cases}
\end{equation}
and $f^h(\bm{x})$ is a smooth blending function that satisfies
\begin{equation}
    f^h(\bm{x}) \jump{\bm{u}}(\bm{x}) = \sum_{I \in n_{nodes}} N_I(\bm{x}) H_{\Gamma}(\bm{x}_I) \jump{\bm{u}}(\bm{x}_I)
\end{equation}
with $n_{nodes}$ denoting the number of nodes in an element.

Enriched finite elements are then constructed by interpolating $\bar{\bm{u}}(\bm{x})$ using the standard (continuous Galerkin) finite element method and interpolating $\jump{\bm{u}}(\bm{x})$ through enrichment. 
The standard finite element interpolation of $\bar{\bm{u}} (\bm{x})$ may be written as two equivalent ways~\cite{cusini2021simulation}: (i) one using scalar shape functions $N_I(\bm{x})$ and vector-valued coefficients $\bm{d}_I$ for the $I$-th node, and (ii) the other using vector-valued shape functions $\bm{N}_i(\bm{x})$ and scalar coefficients $d_i$ for the $i$-th degree of freedom in standard finite elements.
So we can write the standard finite element approximation of $\bar{\bm{u}} (\bm{x})$ as
\begin{equation}
    \bar{\bm{u}} (\bm{x}) \approx \sum_{I \in n_{nodes}} N_I (\bm{x}) \bm{d}_I = \sum_{i \in n_{std}} \bm{N}_i(\bm{x}) d_i,
    \label{eq:interpolation-continuous-u}
\end{equation}
where $n_{std}$ is the number of degrees of freedom in standard finite elements. Next, enrichment of $\jump{\bm{u}}(\bm{x})$ can be written as
\begin{equation}
    \jump{\bm{u}}(\bm{x}) = \sum_{j \in n_{enr}} \bm{\xi}_j(\bm{x}) a_j,
    \label{eq:interpolation-discontinuous-u}
\end{equation}
where $\bm{\xi}_j(\bm{x})$ and $a_j$ are the vector-valued shape functions and scalar coefficients for the $j$-th enriched degree of freedom, respectively, and $n_{enr}$ is the number of enriched degrees of freedom.

Inserting Eqs.~\eqref{eq:interpolation-continuous-u} and~\eqref{eq:interpolation-discontinuous-u} into Eq.~\eqref{eq:displacement-decomposition} gives the enriched finite element approximation of $\bm{u}(\bm{x})$ as
\begin{linenomath}
\begin{align}
    \bm{u} (\bm{x}) 
    &\approx \sum_{i \in n_{std}} \bm{N}_i(\bm{x}) d_i + M_{\Gamma}(\bm{x}) \sum_{j \in n_{enr}} \bm{\xi}_j (\bm{x}) a_j 
    = \sum_{i \in n_{std}} \bm{N}_i(\bm{x}) d_i + \sum_{j \in n_{enr}} \bm{\phi}_j(\bm{x}) a_j,
    \label{eq:displacement-expansion}
\end{align}
\end{linenomath}
where a new function $\bm{\phi}_j(\bm{x})$ is defined as
\begin{equation}
    \bm{\phi}_j(\bm{x}) := H_{\Gamma}(\bm{x}) \bm{\xi}_j(\bm{x}) - \sum_{I \in n_{nodes}} N_I(\bm{x}) H_{\Gamma}(\bm{x}_I) \bm{\xi}_j (\bm{x}_I).
\end{equation}

The foregoing formulation is general for a wide class of enriched finite element methods, and the particular type of enriched method is determined by the specific form of $\bm{\phi}_j$. 
In what follows, we specialize the formulation to XFEM.

\subsection{Extended finite element method}
The XFEM approach enriches global degrees of freedom to interpolate the discontinuous part of the displacement field, $\jump{\bm{u}}(\bm{x})$, in the same way as the continuous part, $\bar{\bm{u}}(\bm{x})$.
This enrichment allows XFEM to describe the complete kinematic modes of the discontinuities, unlike locally enriched AES methods where some kinematic modes are neglected by design.
At the expense of this feature, however, XFEM entails a larger number of enriched degrees of freedom than the AES methods.

Since XFEM interpolates $\jump{\bm{u}}(\bm{x})$ in the same manner as the standard finite element method, $n_{enr} = n_{std}$, $\bm{\xi}_j(\bm{x}) = \bm{N}_j(\bm{x})$, and $\bm{\phi}_j(\bm{x})$ takes the following form
\begin{equation}
    \bm{\phi}_j(\bm{x}) = [H_{\Gamma}(\bm{x}) - H_{\Gamma}(\bm{x}_j)] \bm{N}_{j}(\bm{x}).
\end{equation}
Because XFEM employs the standard Galerkin method, it approximates the variation $\bm{\eta}$ in the same way as the $\bm{u}$, \ie
\begin{equation}
    \bm{\eta} = \bar{\bm{\eta}} + M_{\Gamma} \jump{\bm{\eta}}.
\end{equation}
Note that $\bar{\bm{\eta}}$ and $\jump{\bm{\eta}}$ are two independent variations.
The variation can be further expanded like Eq.~\eqref{eq:displacement-expansion}, and it is omitted for brevity.

Substituting the enriched finite element approximations of $\bm{u}$ and $\bm{\eta}$ into Eq.~\eqref{eq:variational-form}, we arrive at the following two discrete variational equations (in residual form):
\begin{linenomath}
\begin{align}
    \bm{\mathcal{R}}_i^{d} &:= 
    \int_{\Omega \setminus \Gamma} \symgrad \bm{N}_i : \bm{\sigma} \, \dd V 
    % - \int_{\Omega} \bm{N}_i \cdot \bm{f} \, \dd V 
    - \int_{\partial_t \Omega} \bm{N}_i \cdot \bar{\bm{t}} \, \dd A
    = 0, \quad i = 1,2,\ldots, n_{std},
    \label{eq:weak-form-standard} \\
    \bm{\mathcal{R}}_j^{a} &:= 
    \int_{\Omega \setminus \Gamma} \symgrad \bm{\phi}_j : \bm{\sigma} \, \dd V 
    + \int_{\Gamma} \bm{N}_j \cdot \bm{t} \, \dd A 
    % - \int_{\Omega} \bm{\phi}_j \cdot \bm{f} \, \dd V
    = 0, \quad j = 1,2,\ldots, n_{enr}.
    \label{eq:weak-form-enriched}
\end{align}
\end{linenomath}

We use Newton's method to solve Eqs.~\eqref{eq:weak-form-standard} and~\eqref{eq:weak-form-enriched}.
The $k$th Newton iteration proceeds in the following two steps:
\begin{linenomath}
\begin{align}
    \text{solving} &\quad \bm{\mathcal{J}}^k \Delta \bm{X} = - \bm{R}^k, \label{eq:newton-linear-system}\\
    \text{updating} &\quad \bm{X}^{k+1} = \bm{X}^{k} + \Delta \bm{X}.
\end{align}
\end{linenomath}
Here, $\bm{R}$ and $\bm{X}$ are the residual vector and the unknown vector, respectively, which are given by
\begin{linenomath}
\begin{align}
    \bm{R}\left(\bm{X}\right)
    := \begin{Bmatrix}
        \bm{\mathcal{R}}^d\\
        \bm{\mathcal{R}}^a
    \end{Bmatrix}, \quad
    \bm{X}
    := \begin{Bmatrix}
        \bm{d} \\
        \bm{a} 
    \end{Bmatrix},
    \label{eq:residual-unknown}
\end{align}
\end{linenomath}
with $\bm{d}$ and $\bm{a}$ denoting the nodal solution vectors for the standard and enriched degrees of freedom, respectively.
The Jacobian matrix, $\bm{\mathcal{J}}$, is given by
\begin{linenomath}
\begin{align}
    \bm{\mathcal{J}} := \frac{\pd \bm{R}\left(\bm{X}\right)}{\pd \bm{X}} = 
    \begin{bmatrix}
        \pd_{\bm{d}} \bm{\mathcal{R}}^d & \pd_{\bm{a}} \bm{\mathcal{R}}^d \\
        \pd_{\bm{d}} \bm{\mathcal{R}}^a & \pd_{\bm{a}} \bm{\mathcal{R}}^a
    \end{bmatrix}.
    \label{eq:jacobian}
\end{align}
\end{linenomath}
In this work, we use a direct linear solver in each Newton iteration.

In the following, we describe two points that are new in the present work: (i) an integration-point-level algorithm to update the traction vector $\bm{t}$ treated by the barrier method, and (ii) a new surface integration scheme to alleviate traction oscillations that may occur in some problems.

\subsection{Traction update algorithm}
Algorithm~\ref{algo:traction-update} describes how to update the barrier-treated traction vector at an integration point during a Newton iteration.
For notational brevity, all the quantities at the current time steps are written without any subscripts, while the displacement jump at the previous time step is denoted by $\jump{\bm{u}}_{n}$.
The calculation of the contact pressure ($p_{\cn}$) and the tangential traction ($\bm{t}_{\ct}$) was explained in the previous section.
It is noted that the tangential traction is computed without evaluating a yield function, due to the use of the smoothed friction law~\eqref{eq:smoothed-tangential-traction}.
The last step of the algorithm is to compute the tangent operator of the traction with respect to the displacement jump, which is critical to the optimal convergence of Newton's method.
Specifically, the tangent operator can be calculated as:
\begin{linenomath}
\begin{align}
    \mathbb{C}_{\text{interface}} 
    := \dfrac{\partial \bm{t}}{\partial \jump{\bm{u}}}
    &= \dfrac{\partial (-p_N \bm{n})}{\partial \jump{\bm{u}}} 
    + \dfrac{\partial \bm{t}_T}{\partial \jump{\bm{u}}} \nonumber\\
    &= k_{\cn} (\bm{n} \otimes \bm{n})
    - \mu k_{\cn} m(u_T)(\bm{m} \otimes \bm{n}) 
    + \mu p_N \dfrac{\partial m(u_T)}{\partial u_T}(\bm{m} \otimes \bm{m})
    \label{eq:traction-CTO}
\end{align}
\end{linenomath}
where $k_{\cn} := -\pd{p_{\cn}}/\pd{u_{\cn}}$ defined in Eq.~\eqref{eq:barrier-stiffness}.
Notably, when $\|\bm{u}_{T}\|$ is found to be zero (\ie~no slip), one may set $\bm{m}=\bm{0}$ in the succeeding steps.
\begin{algorithm}[h!]
    \setstretch{1.2}
    \caption{Traction update procedure in the barrier-treated enriched finite element method.}
    \begin{algorithmic}[1]
        \Require $\Delta \jump{\bm{u}}$ and $\bm{n}$ at the current step, and $\jump{\bm{u}}_{n}$ at the previous step.
        \State Update the displacement jump in the current step, $\jump{\bm{u}} = \jump{\bm{u}}_{n} + \Delta \jump{\bm{u}}$.
        \State Calculate the normal displacement jump $u_N = \jump{\bm{u}}\cdot\bm{n}$.
        \State Compute the contact pressure, $p_N = -\partial B/\partial u_N$, as in Eq.~\eqref{eq:contact-pressure}.
        \State Calculate the slip displacement vector, $\bm{u}_T  = \jump{\bm{u}} - u_{N}\bm{n}$, the slip magnitude, $u_T = \|\bm{u}_T\|$, and the slip direction vector, $\bm{m}=\bm{u}_T/\|\bm{u}_T\|$.
        \State Compute the tangential traction vector, $\bm{t}_T  = m(u_T)\mu p_N\bm{m}$.
        \State Calculate the traction vector, $\bm{t} = -p_N \bm{n} + \bm{t}_T$.
        \State Compute the tangent operator, $\mathbb{C}_{\text{interface}}$, as in Eq.~\eqref{eq:traction-CTO}.
        \Ensure $\bm{t}$ and $\mathbb{C}_{\text{interface}}$ at the current step.
    \end{algorithmic}
    \label{algo:traction-update}
\end{algorithm}

One can see that the traction update procedure of the barrier method is as simple as that of the classic penalty method -- perhaps even simpler because it does not require one to evaluate a yield function.
We again note that one can easily change a penalty-based code to a barrier-based code, by substituting Algorithm~\ref{algo:traction-update} into an existing traction update algorithm.

\subsection{Averaged surface integration scheme}
As a final task of discretization, we devise a simple method that mitigates spurious oscillations in traction field solutions.
Embedded finite element solutions often exhibit non-physical oscillations in traction fields when the enriched degrees of freedom exert strong constraints on the kinematics.
The origin of such oscillations can be attributed to the inf--sup stability problem common in various types of constrained problems (\eg~\cite{bochev2006stabilization,choo2015stabilized,choo2019stabilized,zhao2020stabilized}).
Several strategies have been proposed to stabilize traction oscillations in embedded finite elements employing the Lagrange multiplier or the penalty method~\cite{ji2004strategies,bechet2009stable,liu2010stabilized}.
While these methods have been found effective, their implementation often requires non-trivial effort.

Here, we present a simple approach to stabilization of traction oscillations in embedded finite elements that employ penalty-like methods including the barrier method.
The approach recasts the weak penalty formulation proposed by Svenning~\cite{svenning2016weak} -- originally developed for a mixed discretization of non-embedded finite elements -- to an XFEM discretization.
In essence, the weak penalty method evaluates the surface integral in the variational equation after projecting $\bm{\tilde{u}}$ to piecewise constant space.
The projection operator is defined as
\begin{equation}
    \mathcal{P} (\tilde{\bm{u}}) := \dfrac{1}{|\Gamma^e|} \int_{\Gamma^e} \tilde{\bm{u}} \, \dd A, \;\; \text{where}\;\;
    |\Gamma^e| := \int_{\Gamma^e} \, \dd A
    \label{eq:projection-opereator-original}
\end{equation}
with $\Gamma^e$ denoting the interface element in non-embedded finite elements.
Importantly, when the displacement field is approximated by linear shape functions, the projection should be performed over two elements across the interface element~\cite{svenning2016weak}.

We now modify the projection operator~\eqref{eq:projection-opereator-original} for XFEM.
As standard, we consider linear approximation of the displacement field, focus on an element cut by a single interface, and introduce two integration points at the interface to evaluate the surface integral in Eq.~\eqref{eq:weak-form-enriched}.
In this case, the cracked element may be viewed as a collection of two displacement elements (where volume integrals are evaluated) and an interface element (where surface integrals are evaluated).
Applying the projection operator~\eqref{eq:projection-opereator-original} to the cracked element then gives a constant value of $\bm{\tilde{u}}$ over the element.
To express this idea mathematically, let $\bm{x}_{Q1}$ and $\bm{x}_{Q2}$ denote the position vectors of the two surface integration points, $A_{Q1}$ and $A_{Q2}$ denote the surface areas calculated at the respective integration points, and $\jump{\bm{u}}_{Q1}$ and $\jump{\bm{u}}_{Q2}$ denote the displacement jump fields at the points.
Then the project operator for XFEM can be written as
\begin{linenomath}
\begin{align}
    \mathcal{P}_{\text{XFEM}}(\tilde{\bm{u}}) 
    &= \dfrac{1}{A_{Q1} + A_{Q2}} [A_{Q1}\tilde{\bm{u}}_{Q1} + A_{Q2} \tilde{\bm{u}}_{Q2} ] \nonumber \\ 
    &= \dfrac{1}{A_{Q1} + A_{Q2}} [A_{Q1}\sum_{I \in n_{nodes}} N_I (\bm{x}_{Q1}) \tilde{\bm{u}}_I + A_{Q2}\sum_{I \in n_{nodes}} N_I (\bm{x}_{Q2}) \tilde{\bm{u}}_I] \nonumber \\
    &= \dfrac{1}{2} \sum_{I \in n_{nodes}} [ N_I (\bm{x}_{Q1}) +  N_I (\bm{x}_{Q2}) ] \tilde{\bm{u}}_I.
\end{align}
\end{linenomath}
In words, the operator averages the two displacement-jump fields calculated from the two integration points.
So we shall refer to this method as an averaged (surface) integration scheme.

To summarize, we stabilize traction oscillations using the averaged integration scheme
\begin{equation}
    \jump{\bm{u}}_{Q1} = \jump{\bm{u}}_{Q2} = \mathcal{P}_{\text{XFEM}}(\tilde{\bm{u}}),
\end{equation}
which is an XFEM version of the weak penalty formulation proposed by Svenning~\cite{svenning2016weak}. 
The scheme uses a single value of the displacement jump to evaluate the surface integral in Eq.~\eqref{eq:weak-form-enriched} within a cracked element.
As will be demonstrated later, the averaged integration scheme can effectively suppress traction oscillations in XFEM.
The key advantage of this method over existing methods is that it can be implemented with trivial effort.

It would be worthwhile to note two points regarding the averaged integration scheme.
First, the scheme is applicable to general penalty-like methods that share the same mathematical structure.
Indeed, it works well for the barrier method because the barrier method has the same mathematical structure as the penalty method, as discussed earlier.
Second, when the standard two-point integration scheme does not show traction oscillations, the averaged scheme would be sub-optimal in terms of accuracy.
So we propose the averaged scheme as a new option for problems that the standard scheme gives severe traction oscillations, rather than an absolute alternative to the standard scheme for general contact problems.
Later in numerical examples, we use the standard integration scheme by default, and switch to the averaged scheme if severe oscillations are observed.

% SECTION
% ==============================================================================
\section{Parameterization}
\label{sec:parameterization}

In this section, we tailor the parameters of the barrier method to embedded interfaces, such that the method can be used without parameter tuning.
The motivation is twofold: (i) application of the barrier method to embedded interfaces requires a new parameter that represents the initial ``gap'' of a closed interface, and (ii) the existing way to determine the barrier stiffness parameter ($\kappa$) cannot be applied to embedded interfaces under quasi-static conditions.
In what follows, we develop a way to constrain these two parameters for a given interface embedded in finite elements. Then the barrier method has only one free parameter, $\hat{d}$, which controls the solution accuracy ($\hat{s}=\hat{d}$ is assumed).
The value of $\hat{d}$ is then recommended as a function of the problem domain size, which gives sufficiently accurate numerical solutions in most cases.

To begin, let us denote by $d_0$ the initial gap of a closed interface.
By definition, $d_0$ should range in between $0$ and $\hat{d}$, but it is unclear how to set this value within this range.
To answer this question, we recall that the gap distance determines the contact pressure and the normal stiffness of a barrier-treated interface.
So we can obtain the initial values of the contact pressure $p_{\cn}$ and the normal stiffness $k_{\cn}$ by inserting $d_0$ into Eq.~\eqref{eq:contact-pressure} and~\eqref{eq:barrier-stiffness}, respectively.
This operation gives the initial contact pressure, $p_{\cn,0}$ and the initial normal stiffness, $k_{\cn,0}$, as follows:
\begin{linenomath}
\begin{align}
    p_{\cn,0} &= \kappa(d_{0} - \hat{d})\left[2\ln\left(\dfrac{d_{0}}{\hat{d}}\right) - \dfrac{\hat{d}}{d_{0}} + 1\right], \label{eq:initial-contact-pressure} \\
    k_{\cn,0} &= -2\kappa\ln\left(\dfrac{d_{0}}{\hat{d}}\right) - \kappa\dfrac{(d_{0}-\hat{d})(3d_{0} + \hat{d})}{d_{0}^{2}}. \label{eq:initial-barrier-stiffness}
\end{align}
\end{linenomath}
It can be seen that for a given value of $\hat{d}$, $d_0$ and $\kappa$ together determine $p_{\cn,0}$ and $k_{\cn,0}$. 
Thus we seek to constrain the values of $d_0$ and $\kappa$ by prescribing the values of $p_{\cn,0}$ and $k_{\cn,0}$.

We first prescribe the value of $p_{\cn,0}$ based on a physics-based argument that the optimal value of $p_{\cn,0}$ is the contact pressure acting on the interface if it was perfectly closed in the initial condition.
This optimal value, which we denote by $(p_{\cn,0})_{\text{optimal}}$, can be calculated or approximated based on our knowledge of the geometry and boundary conditions of the problem at hand.
So we postulate that $(p_{\cn,0})_{\text{optimal}}$ is given and assign $p_{\cn,0}=(p_{\cn,0})_{\text{optimal}}$, in order to let the barrier-approximated interface behave like a perfectly closed interface as much as possible.

\begin{figure}[h!]
    \centering
    \includegraphics[width=0.5\textwidth]{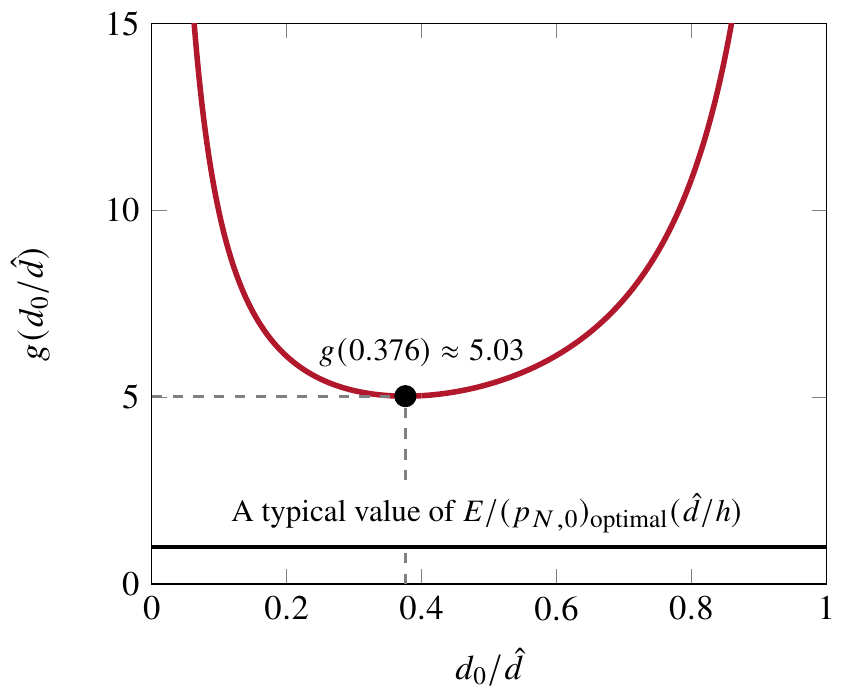}
    \caption{Variation of $g(d_0/\hat{d})$ with $d_0/\hat{d}\in[0,1]$.}
    \label{fig:bowl-shaped-function}
\end{figure}

Next, we prescribe the value of $k_{\cn,0}$ such that it makes the condition of the Jacobian matrix~\eqref{eq:jacobian} as good as possible.
When the bulk material has a constant stiffness, the optimal value of $k_{\cn,0}$ is approximated to be $E/h$, where $E$ is Young's modulus and $h$ is the element size.
See \ref{appendix:A} for the details of the approximation.
We thus seek to find the value of $d_{0}$ that makes $k_{\cn,0}=E/h$ when $(p_{\cn,0})_{\text{optimal}}$ is given. 
To this end, we rearrange Eq.~\eqref{eq:initial-contact-pressure} in terms of $\kappa$, insert the expression for $\kappa$ into Eq.~\eqref{eq:initial-barrier-stiffness}, equate $p_{\cn,0}$ to $(p_{\cn,0})_{\text{optimal}}$, and let $k_{\cn,0}=E/h$.
This gives
\begin{equation}
    g(d_0/\hat{d}) = \dfrac{E}{(p_{\cn,0})_{\text{optimal}}}\dfrac{\hat{d}}{h},
    \label{eq:d0_optimal}
\end{equation}
where
\begin{equation}
    g(d_0/\hat{d}) := \dfrac{-2 \mathrm{ln}(d_0/\hat{d}) - [1-(d_0/\hat{d})^{-1}][3+(d_0/\hat{d})^{-1}]}{(d_0/\hat{d} - 1) \left[ 2 \mathrm{ln} (d_0/\hat{d}) - (d_0/\hat{d})^{-1} + 1 \right]}.
\end{equation}
One can try to solve Eq.~\eqref{eq:d0_optimal} to find the optimal value of $d_0/\hat{d}$.
Unfortunately, however, Eq.~\eqref{eq:d0_optimal} is usually not solvable, because the minimum value of the left hand side, $g(d_0/\hat{d})$, is usually greater than the value of the right hand side, $E/(p_{\cn,0})_{\text{optimal}}(\hat{d}/{h})$.
To illustrate this point, in Fig.~\ref{fig:bowl-shaped-function} we plot the values of $g(d_0/\hat{d})$ in the range of $0\leq d_0/\hat{d} \leq 1$.
As can be seen, $g(d_0/\hat{d})$ is a bowl-shaped function and attains its minimum value of $\sim 5.03$ when $d_0/\hat{d}=0.376$.
However, $E/(p_{\cn,0})_{\text{optimal}}(\hat{d}/{h})$ is typically lower than $\sim 5.03$ because $\hat{d}\ll h$ in practical analysis settings. 
So Eq.~\eqref{eq:d0_optimal} does not have a solution in almost all cases.
Still, however, we can see that $d_0/\hat{d}=0.376$ makes the value of $g(d_0/\hat{d})$ closest to the optimal value, which means that it will make the matrix condition as good as possible under the given physical and numerical parameters.

Based on the foregoing statement, we set $d_0$ as
\begin{equation}
    d_{0} = 0.376\,\hat{d}.
\end{equation}
Also, rearranging Eq.~\eqref{eq:initial-contact-pressure} and substituting $(p_{\cn,0})_{\rm optimal}$ into $p_{\cn}$, we can write $\kappa$ as a function of $d_0$ and then insert $d_{0} = 0.376\,\hat{d}$ as follows:
\begin{linenomath}
\begin{align}
    \kappa &= \dfrac{(p_{\cn,0})_{\rm optimal}}{(d_{0} - \hat{d})\left[2\ln\left(\dfrac{d_{0}}{\hat{d}}\right) - \dfrac{\hat{d}}{d_{0}} + 1\right]} \label{eq:kappa} \\
    &= \dfrac{(p_{\cn,0})_{\rm optimal}}{2.256\,\hat{d}}.
\end{align}
\end{linenomath}
The upshot of the above parameterization is that now $d_0$ are $\kappa$ constrained to $\hat{d}$.
Then if we take $\hat{d}=\hat{s}$, $\hat{d}$ is the only free parameter of the method.

While the value of $\hat{d}$ can technically be assigned with any positive value according to the desired accuracy of the numerical solution, we also recommend $\hat{d}=10^{-4}L$ where $L$ is the size of the problem domain.
This recommendation is based on our own experience with a variety of numerical examples, as well as on the values of $\hat{d}$ used for multi-body contact problems~\cite{li2020incremental}.
In Li \etal~\cite{li2020incremental}, $\hat{d}=10^{-3}L$ is used for most of their multi-body contact problems.
Yet we have experienced that $\hat{d}=10^{-3}L$ gives slight errors under a large compression, presumably because $d_0$ is not sufficiently small. 
According to our experience, however, $\hat{d}=10^{-4}L$ always give sufficiently accurate solutions. 
So we recommend $\hat{d}=10^{-4}L$ as a rule of thumb for embedded interface problems.
\revised{Note also that a lower value of $\hat{d}$ can be well used if desired (\eg~when the domain size is very large like a geologic scale.)}

Table~\ref{tab:barrier-parameters} summarizes the parameters of the barrier method and their recommended values for an embedded interface problem.
Unless otherwise specified, we use these parameters in the numerical examples in the next section.
\begin{table}[h!]
    \centering
    \begin{tabular}{llllr}
    \toprule
    Parameter & Symbol & Recommended Value \\
    \midrule
    Barrier thickness & $\hat{d}$ & $10^{-4}L$\\
    Maximum microslip & $\hat{s}$ & $\hat{d}$ \\
    Initial gap distance & $d_{0}$ & $0.376\,\hat{d}$ \\
    Barrier stiffness & $\kappa$  &$(p_{\cn,0})_{\rm optimal}/(2.256\,\hat{d})$ \\
    \bottomrule
    \end{tabular}
    \caption{Parameters of the barrier method and their recommended values for an embedded interface problem. The user can first set the value of $\hat{d}$ according to the domain size $L$, and then use the value of $\hat{d}$ to determine the other parameters. It is noted that $(p_{\cn,0})_{\rm optimal}$ is calculated based on the geometry and boundary conditions of the problem at hand, and that $(p_{\cn,0})_{\rm optimal}$ does not need to be calculated accurately (see Remark~\ref{rem:p_0_optimal}).}
    \label{tab:barrier-parameters}
\end{table}
\smallskip

\begin{remark}
    In practice, it suffices to use an approximate value of $(p_{\cn,0})_{\rm optimal}$ to calculate $\kappa$ because $\kappa$ controls the condition of the Jacobian matrix~\eqref{eq:jacobian} but not the solution accuracy.
    According to our experience, even when the value of $(p_{\cn,0})_{\rm optimal}$ is one or two orders of magnitude higher than that of the ``true'' value, the matrix is sufficiently well-conditioned.
    We have observed, however, if the assigned value for $(p_{\cn,0})_{\rm optimal}$ (and hence $\kappa$) is too low, sometimes the Newton iteration does not converge at early steps.
    It is therefore recommended that one assigns the highest possible value for $(p_{\cn,0})_{\rm optimal}$ to calculate the value of $\kappa$.
    \label{rem:p_0_optimal}
\end{remark}
\smallskip

\begin{remark}
    In a rare case where Eq.~\eqref{eq:d0_optimal} is solvable, one may determine $d_0$ by solving the equation and take the smaller solution. (The equation usually has two solutions, see Fig.~\ref{fig:bowl-shaped-function}.) Then the solution can be used to calculate $\kappa$ through Eq.~\eqref{eq:kappa}.
\end{remark}

\smallskip
\begin{remark}
    The foregoing parameterization procedure has been developed assuming that the user wants to follow the classical model formulation in which asperity effects are neglected. 
    However, if the proposed method is used for incorporating asperity effects (see Remarks~\ref{rem:impossible-zero-gap} and~\ref{rem:microslip}), one can also assign the values of $\hat{d}$ and $\hat{s}$ freely as desired.
    Still, $d_0$ and $\kappa$ may be determined as recommended in Table~\ref{tab:barrier-parameters} to avoid ill-conditioning.
\end{remark}

% SECTION
% ==============================================================================
\section{Numerical examples}
\label{sec:numerical-examples}
This section has two purposes: 
(i) to verify the barrier method for embedded interfaces,
and (ii) to investigate the performance of the barrier method under diverse conditions.
For these purposes, we apply the barrier method to numerical examples with varied levels of complexity, from a simple horizontal crack to a circular interface between two materials.
In doing so, we compare the accuracy and robustness of the barrier method with those of the penalty method which has more or less the same computational cost.
The penalty method is incorporated into XFEM following the formulation of Liu and Borja~\cite{liu2008contact} in which a normal penalty parameter, $\alpha_{\cn}$, and a tangential penalty parameter, $\alpha_{\ct}$, are used as normal and tangential spring constants, respectively. 
We note that such a penalty method has been commonly used in XFEM and compared with other contact algorithms, see, \eg~\cite{khoei2007enriched,liu2010finite,annavarapu2014nitsche}.
It is however noted that whenever the contact pressure is positive, the penalty method unavoidably permits inter-penetration ($u_{\cn}<0$) because it calculates the contact pressure as $p_{\cn}=\alpha_{\cn}(-u_{\cn})$.

The numerical examples in this section are prepared using an in-house finite element code \texttt{Geocentric}, which is built on the \texttt{deal.II} finite element library~\cite{arndt2021deal}.
Verification of the code for a variety of non-interface problems can be found from the literature (\eg~\cite{choo2018enriched,choo2018large}).
In the following, we assume plane-strain conditions and  use quadrilateral elements with linear shape functions.

\subsection{Horizontal crack under compression and shear}
\begin{figure}[h!]
    \centering
    \includegraphics[width=0.55\textwidth]{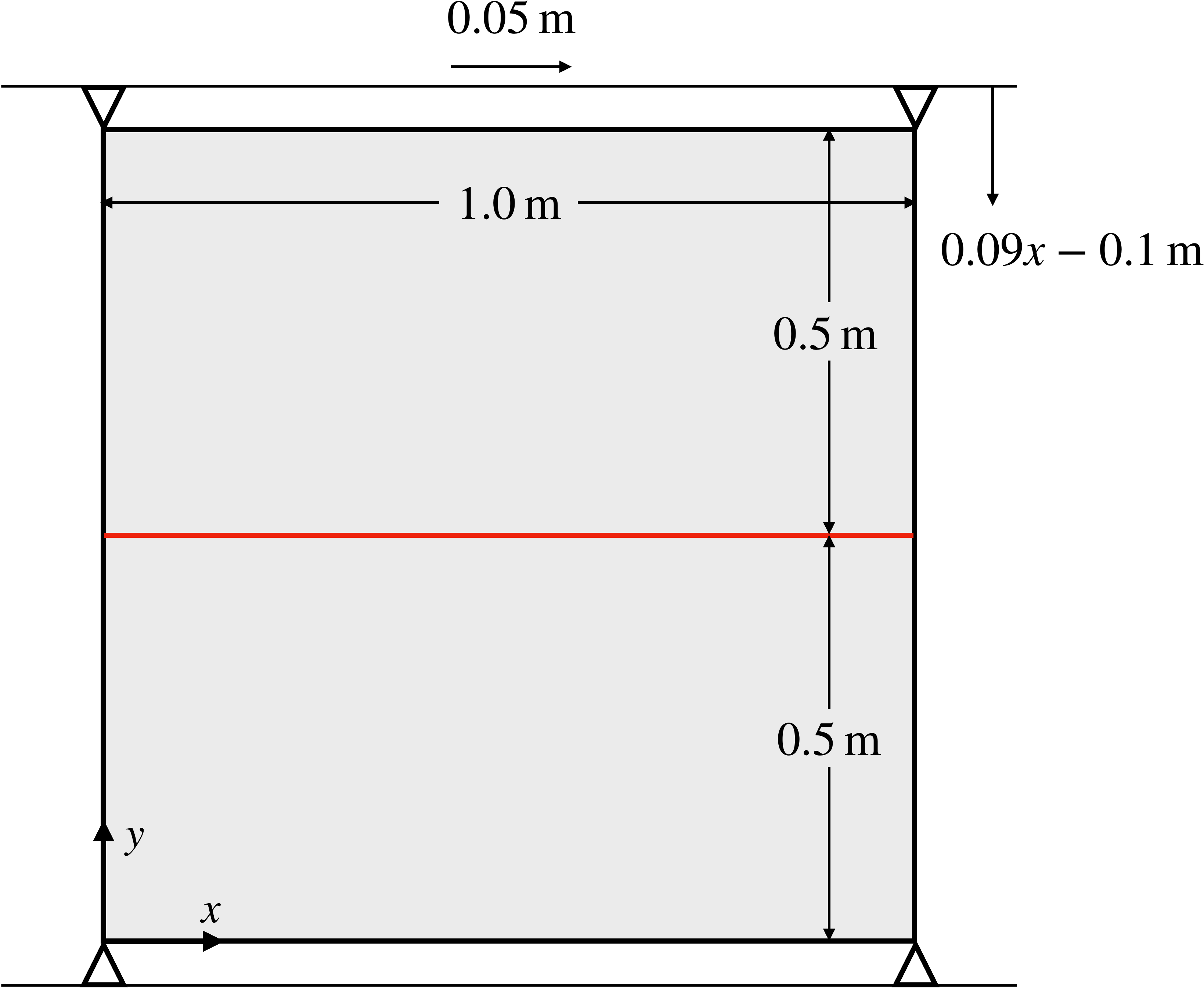}
    \caption{Horizontal crack under compression and shear: problem geometry and boundary conditions.}
    \label{fig:horizontal-setup}
\end{figure}

Our first example applies the barrier method to the problem of a horizontal crack under compression and shear, which was simulated by Liu and Borja~\cite{liu2008contact} and Annavarapu \etal~\cite{annavarapu2014nitsche} using XFEM.
Figure~\ref{fig:horizontal-setup} depicts the geometry and boundary conditions of the problem.
The friction coefficient of the crack is $\mu=0.3$, and the Young's modulus and Poisson's ratio of the bulk material are $E=10$ GPa and $\nu=0.3$, respectively.
We calculate the barrier method parameters as suggested in Table~\ref{tab:barrier-parameters}, with $(p_{N,0})_{\text{optimal}} = 0.55$ GPa which is calculated assuming that the top boundary's displacement is uniform and equal to the average of the actual displacement.
To conduct a mesh refinement study with the embedded crack, we discretize the domain using three different meshes: Mesh 1 ($h=1/11$ m), Mesh 2 ($h=1/25$ m), and Mesh 3 ($h=1/51$ m). 
Each mesh is comprised of square elements of uniform size.

\begin{figure}[h!]
    \centering
    \subfloat[Normal displacement jump]{\includegraphics[width=0.45\textwidth]{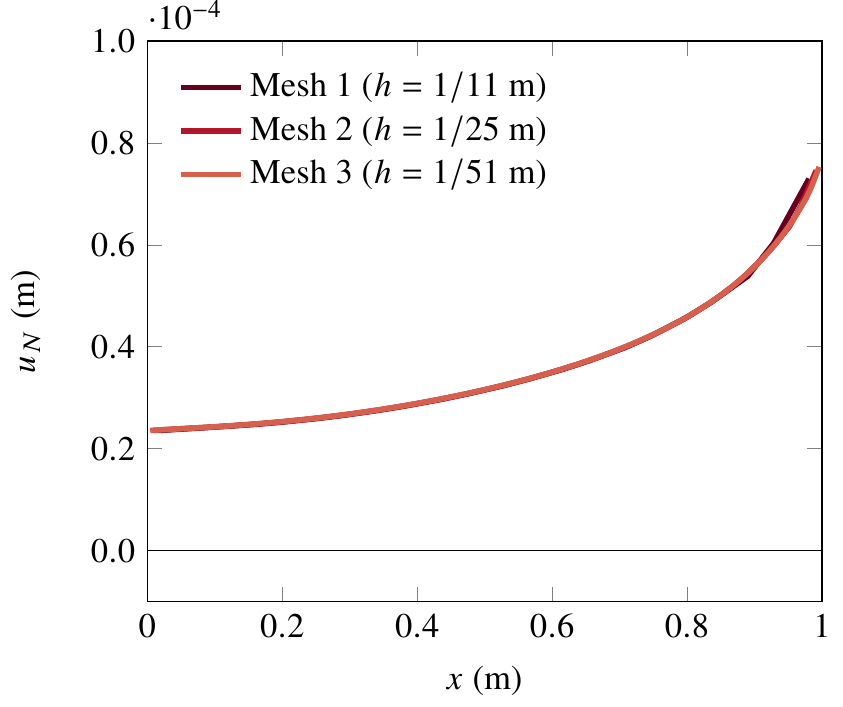}}$\quad$
    \subfloat[Slip displacement]{\includegraphics[width=0.45\textwidth]{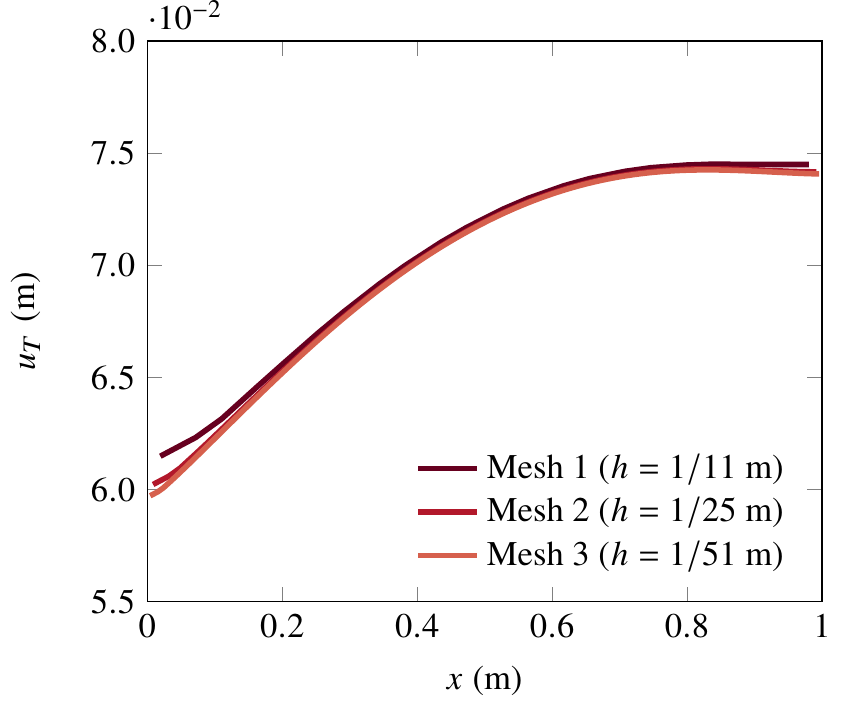}}\\
    \subfloat[Contact pressure]{\includegraphics[width=0.45\textwidth]{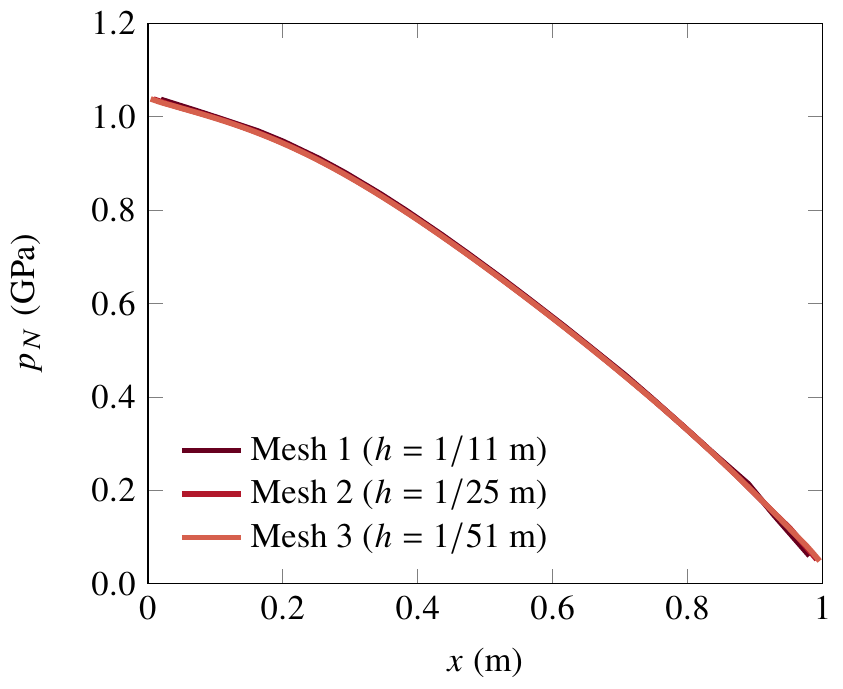}}$\quad$
    \subfloat[Tangential stress]{\includegraphics[width=0.45\textwidth]{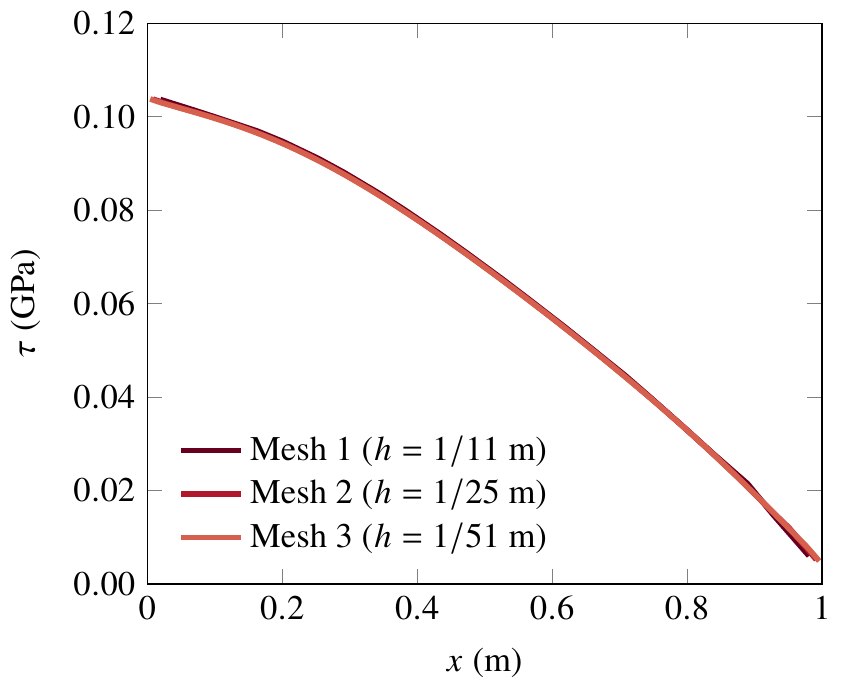}}
    \caption{Horizontal crack under compression and shear: mesh refinement study.}
    \label{fig:horizontal-mesh}
\end{figure}

Figure~\ref{fig:horizontal-mesh} presents the numerical solutions obtained by the three different meshes, in terms of the normal displacement jump ($u_{\cn}$), slip displacement ($u_{\ct}$), contact pressure ($p_{\cn}$), and the tangential stress ($\tau$).
It can be seen that all these quantities converge very well upon mesh refinement.
One may also find that the converged solutions agree well with those obtained by Annavarapu \etal~\cite{annavarapu2014nitsche} using the weighted Nitsche method and the penalty method.

\begin{figure}[h!]
    \centering
    \subfloat[Barrier]{\includegraphics[width=0.45\textwidth]{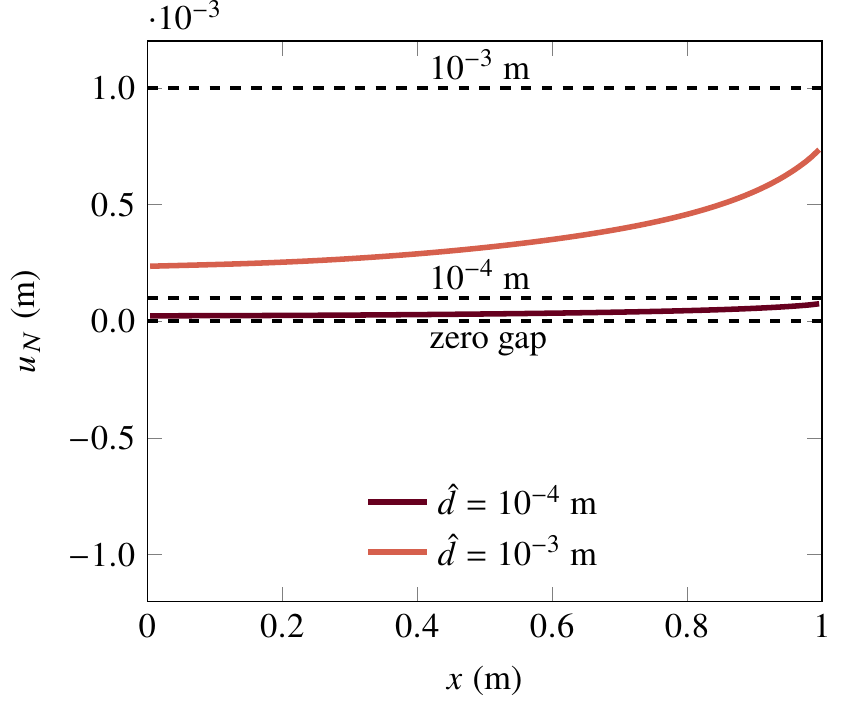}}$\quad$
    \subfloat[Penalty]{\includegraphics[width=0.45\textwidth]{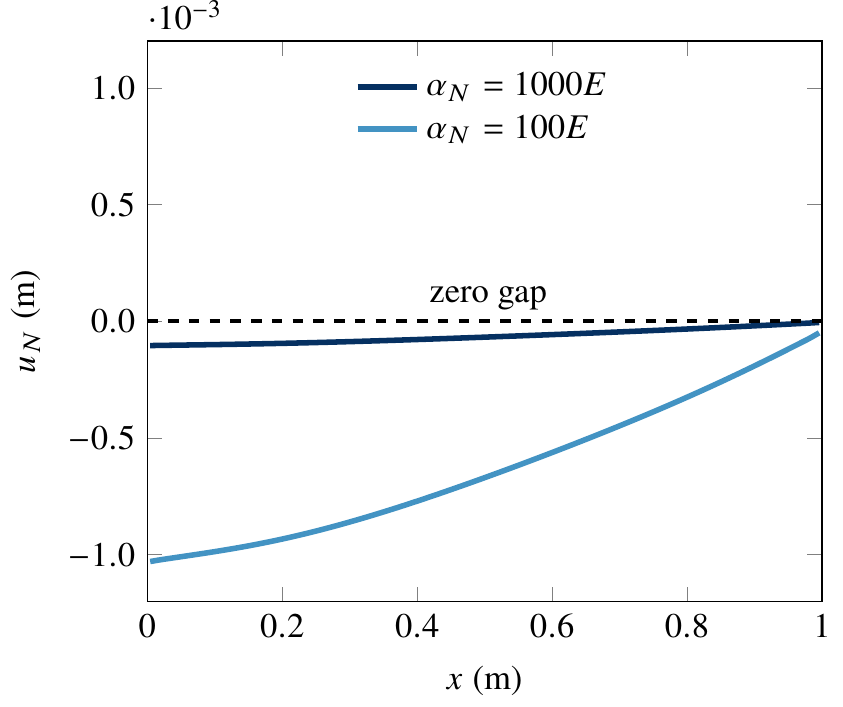}}
    \caption{Horizontal crack under compression and shear: comparison between the barrier and penalty solutions to the normal displacement jump.}
    \label{fig:horizontal-barrier-penalty-uN}
\end{figure}

In Fig.~\ref{fig:horizontal-barrier-penalty-uN}, we demonstrate how the normal displacement jump -- a direct measure of numerical error -- of the barrier method changes by its parameter, in comparison with the penalty method.
When the barrier thickness parameter $\hat{d}$ decreases from $10^{-4}$ m to $10^{-3}$ m, the gap distance in the barrier solution increases.
Likewise, when the normal penalty parameter $\alpha_{\cn}$ decreases by an order or magnitude, the amount of interpenetration in the penalty solution increases.
An important difference, however, exists between the barrier and penalty solutions: the maximum error in the barrier solution is guaranteed to be less than $\hat{d}$, whereas that in the penalty solution cannot be controlled directly.  

\begin{figure}[h!]
    \centering
    \subfloat[Barrier]{\includegraphics[width=0.45\textwidth]{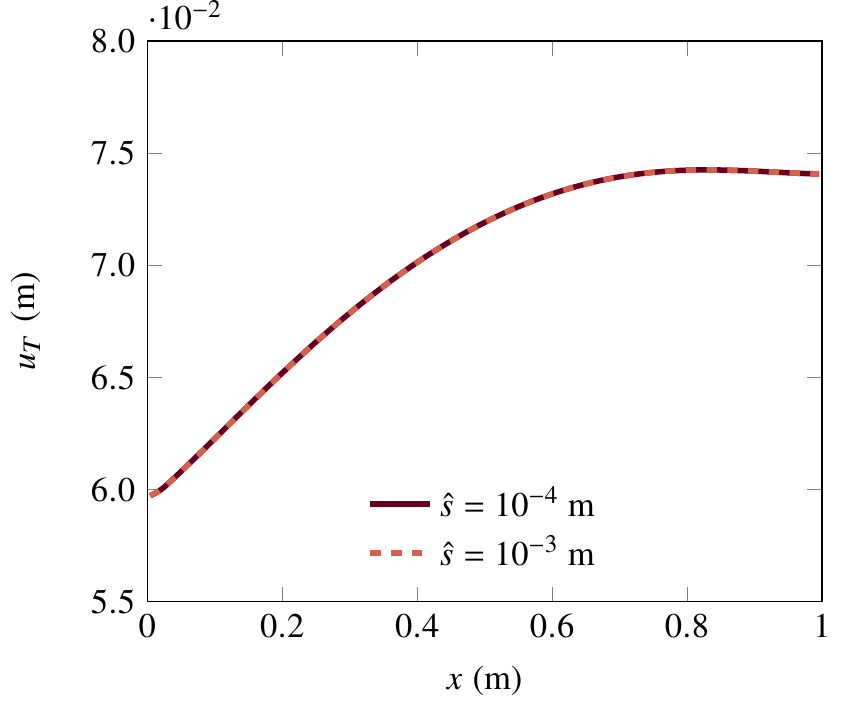}}$\quad$
    \subfloat[Penalty]{\includegraphics[width=0.45\textwidth]{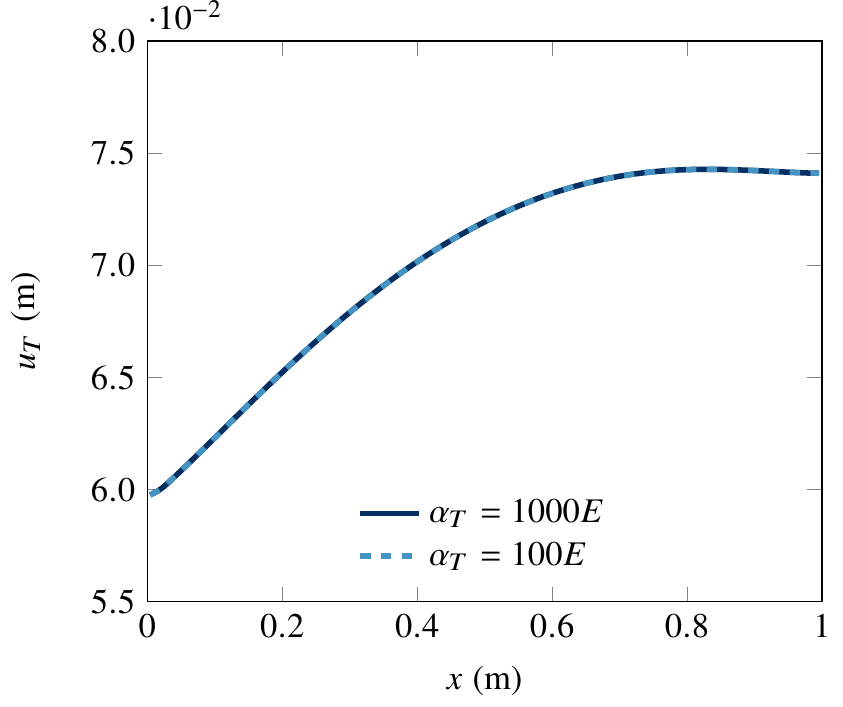}}
    \caption{Horizontal crack under compression and shear: comparison between the barrier and penalty solutions to the slip displacement.}
    \label{fig:horizontal-barrier-penalty-uT}
\end{figure}

In Fig.~\ref{fig:horizontal-barrier-penalty-uT} we further examine how the slip displacement solutions of the barrier and penalty methods are sensitive to their parameters.
As can be seen, the slip displacement solutions are virtually insensitive to the parameters.
The reason would be that the magnitudes of slip displacements are much larger than those of numerical errors.
While not presented for brevity, we have found that the slip displacement solution of the barrier method exhibits noticeable errors only when $\hat{s}$ is larger than the maximum amount of slip displacement.

\subsection{Inclined crack under compression}
\begin{figure}[htbp]
    \centering
    \includegraphics[width=0.45\textwidth]{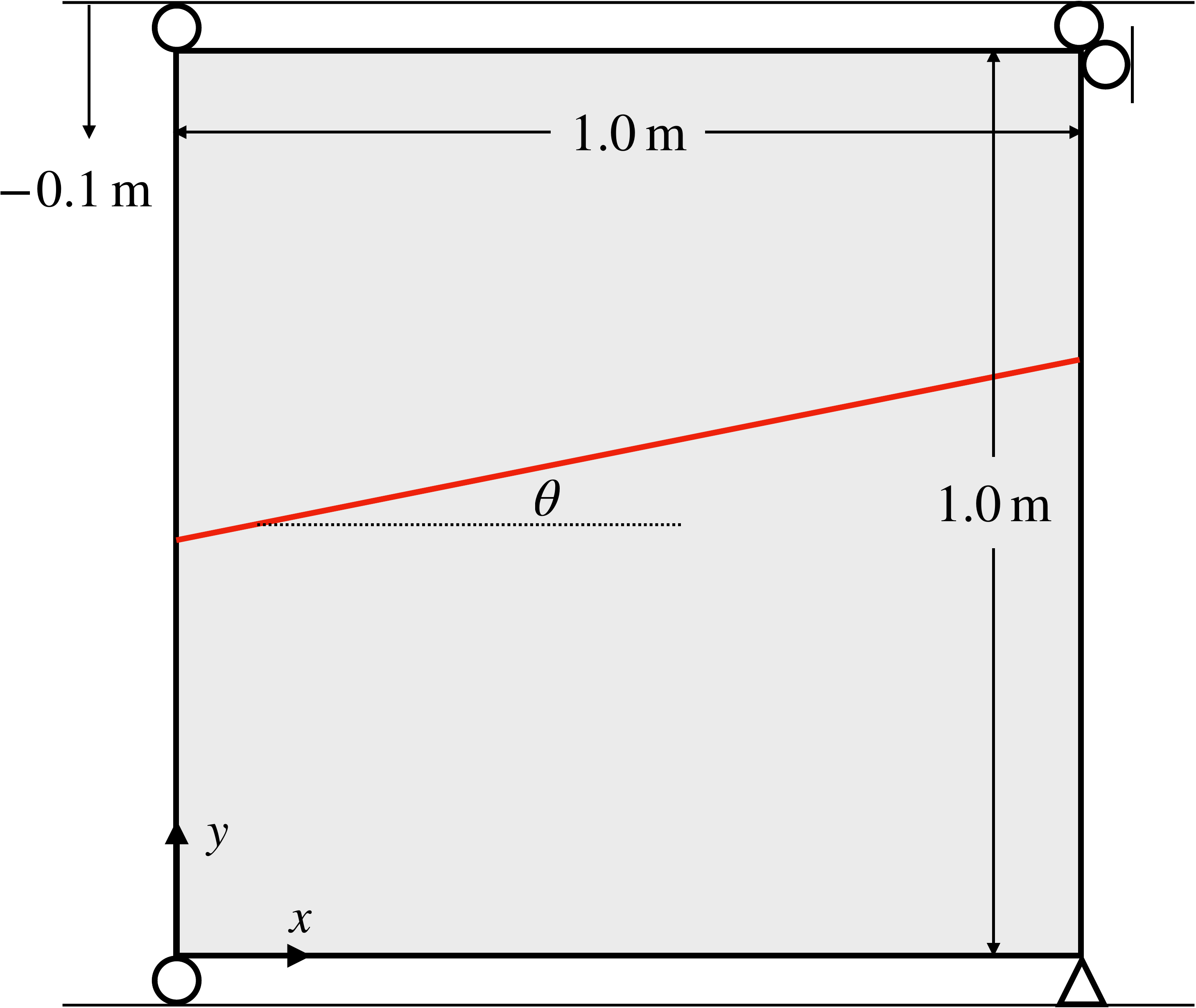}
    \caption{Inclined crack under compression: problem geometry and boundary conditions.}
    \label{fig:inclined-crack-setup}
\end{figure}

In our second example, we examine the (in)sensitivity of the barrier method to the value of $(p_{\cn,0})_{\text{optimal}}$, which is used to calculate the parameter $\kappa$ as suggested in Table~\ref{tab:barrier-parameters}.
The motivation is that unlike the previous problem, it is sometimes not straightforward to calculate $(p_{\cn,0})_{\text{optimal}}$.
As an example, we adopt the problem of compression of a square domain with an inclined crack, which was first presented in Dolbow \etal~\cite{dolbow2001extended} and later used by other studies on embedded frictional interfaces (\eg~\cite{annavarapu2014nitsche,fei2020phase-a}).
Figure~\ref{fig:inclined-crack-setup} illustrates the geometry and boundary conditions of the problem.  
The interface is inclined with an angle $\theta = \text{tan}^{-1}(0.2)$. 
The friction coefficient of the interface is set as $\mu=0.19$ so that it undergoes slip under compression. Following the original problem, we assign the Young's modulus and Poisson's ratio of the material as $E=1$ GPa and $\nu=0.3$, respectively.
We discretize the domain by $160$ by $160$ square elements and apply the compression through 10 steps.

Because the value of $(p_{\cn,0})_{\text{optimal}}$ in this problem is not constant along the interface and not straightforward to calculate, we consider two options.
The first option is to use the upper bound on the initial contact pressure that would have exerted on the interface if it was horizontal.
Specifically, the upper bound value is calculated to be $(p_{N,0})_{\rm horizontal}=10$ MPa under the compressive displacement at the first load step.
As the second option, we use $0.5(p_{N,0})_{\rm horizontal}$ to approximate $(p_{\cn,0})_{\text{optimal}}$.
These two options will lead to two different values of $\kappa$.
Note, however, that all other parameters are unaffected and set according to Table~\ref{tab:barrier-parameters}. 

\begin{figure}[h!]
    \centering
    \includegraphics[width=0.8\textwidth]{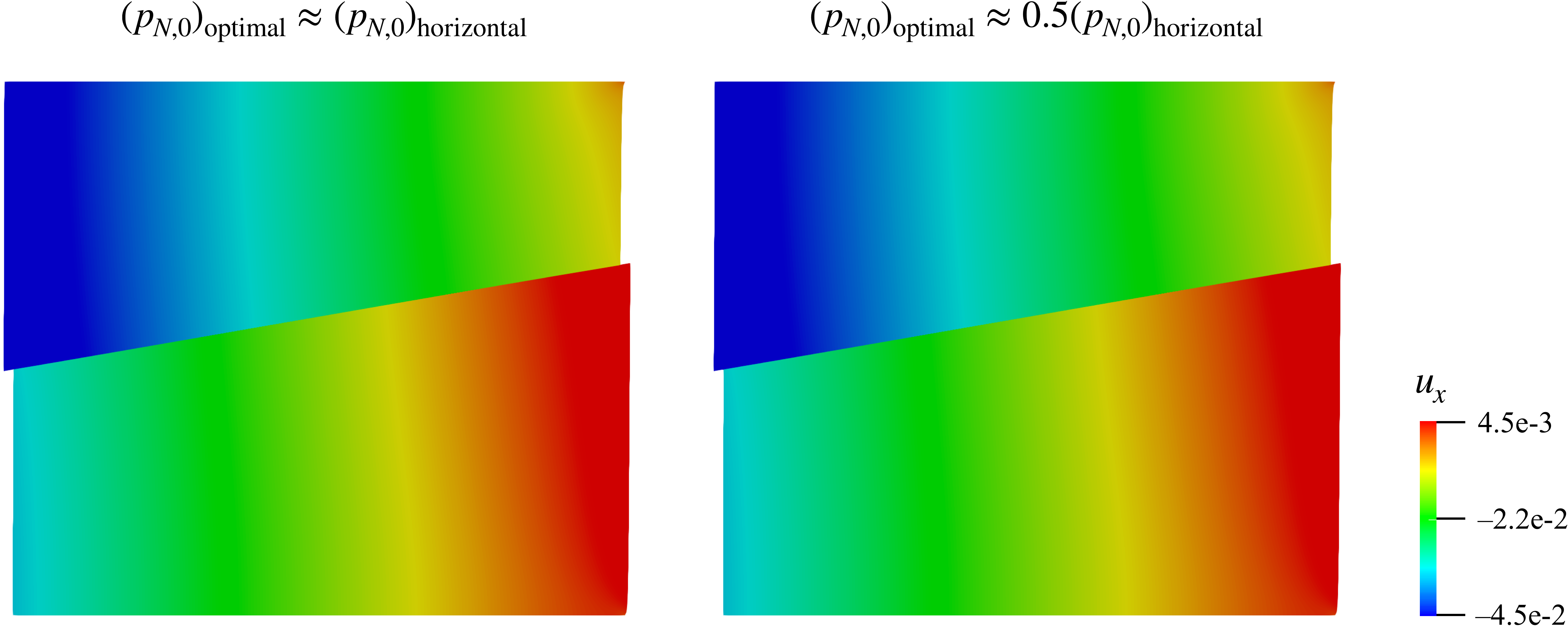}
    \caption{Inclined crack under compression: the $x$-displacement fields obtained with two different approximations to $(p_{N,0})_{\rm optimal}$ in calculating $\kappa$. The domains are deformed by the displacement fields. Color bar in meters.}
    \label{fig:inclined-crack-x-disp}
\end{figure}

\begin{figure}[h!]
    \centering
    \subfloat[Contact pressure]{\includegraphics[width=0.45\textwidth]{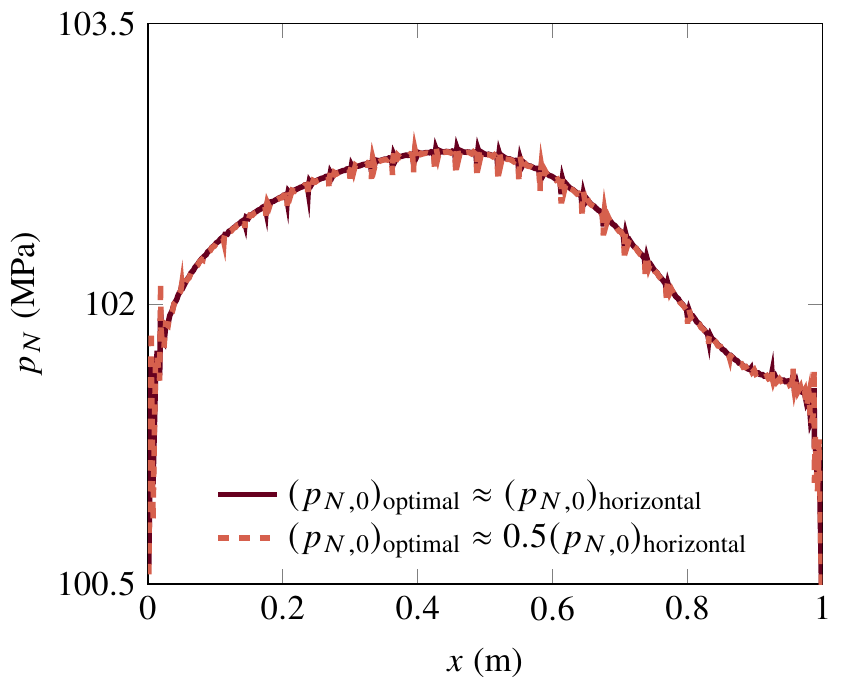}}$\quad$
    \subfloat[Tangential stress]{\includegraphics[width=0.45\textwidth]{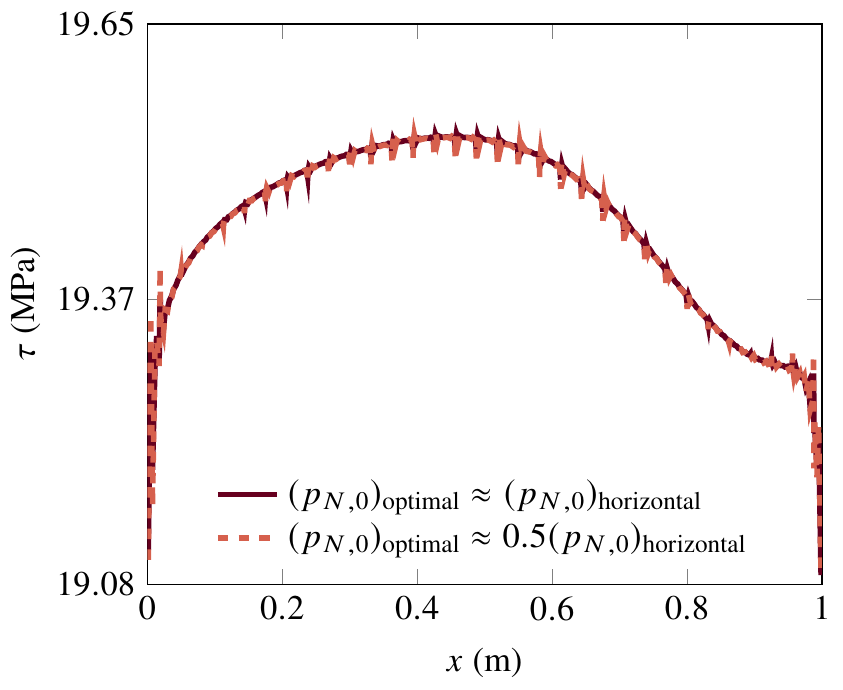}}
    \caption{Inclined crack under compression: (a) contact pressures and (b) tangential stresses obtained with two different approximations to $(p_{N,0})_{\rm optimal}$ in calculating $\kappa$.}
    \label{fig:inclined-crack-pN-tau}
\end{figure}

In Figs.~\ref{fig:inclined-crack-x-disp} and~\ref{fig:inclined-crack-pN-tau}, we compare the numerical solutions obtained with the two different approximations to $(p_{N,0})_{\rm optimal}$, in terms of the displacement field and the interface traction (contact pressure and tangential stress), respectively.
One can see that the two solutions are practically identical.
While not presented for brevity, we have confirmed that these solutions are virtually the same as the penalty solutions to the same problem.
Notably, we have also observed that the solutions are unaffected when $(p_{N,0})_{\rm optimal}$ was approximated with values even much higher than $(p_{N,0})_{\rm horizontal}$.
This is not surprising because $(p_{N,0})_{\rm optimal}$ and thus $\kappa$ controls the matrix condition only -- not the solution accuracy -- as demonstrated in Figs.~\ref{fig:inclined-crack-x-disp} and~\ref{fig:inclined-crack-pN-tau}.
Yet we have also observed that when $(p_{N,0})_{\rm optimal}$ is approximated to be lower than $0.5(p_{N,0})_{\rm horizontal}$, the Newton iteration did not converge, presumably because the barrier is too ``soft'' in some places.
Therefore, as described in Remark~\ref{rem:p_0_optimal}, we recommend using an upper bound value to approximate $(p_{N,0})_{\rm optimal}$ in calculating the value of $\kappa$. 

Before moving to the next example, we note that the traction solutions in Fig.~\ref{fig:inclined-crack-pN-tau} show slight oscillations like those obtained by penalty and other methods (see Annavarapu \etal~\cite{annavarapu2014nitsche} for example).
As the oscillations are not so significant, we have not used the averaged surface integration scheme for this example.
In the last example of this section, we will deal with a more challenging problem that shows much severe oscillations and apply the averaged integration scheme to demonstrate its performance.

\subsection{Sliding between two blocks with stiffness contrast}
\begin{figure}[h!]
    \centering
    \includegraphics[width=0.55\textwidth]{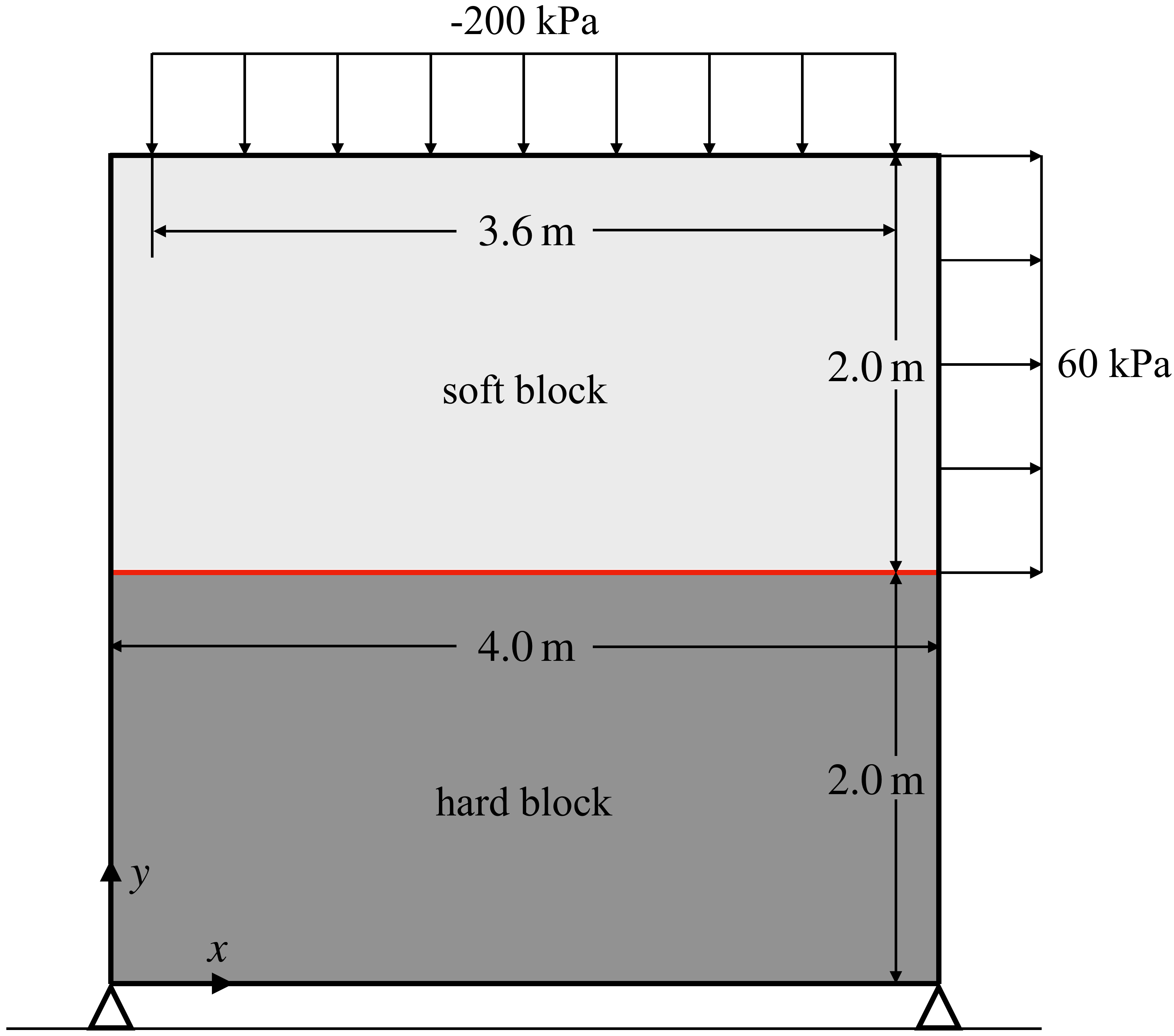}
    \caption{Sliding between two blocks with stiffness contrast: problem geometry and boundary conditions.}
    \label{fig:stiffness-contrast-setup}
\end{figure}

Our third example investigates the performance of the smoothed friction law for handling challenging sliding problems. To this end, we simulate sliding between two blocks with stiffness contrast, extending the problem of an elastic block on a rigid foundation that was designed originally by Oden and Pires~\cite{oden1984algorithms} and later used in several other studies (\eg~\cite{simo1992augmented,annavarapu2014nitsche,fei2020phase-a}).
As the original problem is not an embedded interface problem, we modify the problem as in Fig.~\ref{fig:stiffness-contrast-setup}, similar to the modification in Annavarapu \etal~\cite{annavarapu2014nitsche}. 
Fixing the Young's modulus of the soft block to be 1000 kPa (the same as the original problem~\cite{oden1984algorithms}), we consider two cases of stiffness contrast between the two blocks: $E_{\text{hard}}/E_{\text{soft}}=10^{1}$ and $E_{\text{hard}}/E_{\text{soft}}=10^{7}$, where $E_{\text{hard}}$ and $E_{\text{soft}}$ denote the Young's moduli of the hard and soft blocks, respectively.
The latter case -- higher stiffness contrast -- is intended to emulate the original problem where the foundation is rigid. 
The Poisson's ratios of the two blocks are both $\nu=0.3$.
The friction coefficient of the interface is $\mu=0.5$.
It is noted that under this problem setup, stick and slip conditions are mixed along the interface. 
The barrier stiffness parameter $\kappa$ is calculated setting $(p_{N,0})_{\text{optimal}}$ equal to the magnitude of the compressive traction, 200 kPa.
We discretize the domain by $40$ by $41$ elements, embedding the interface inside elements.

The main reason for selecting this problem is that the penalty method is known to not work well for this problem.
Annavarapu \etal~\cite{annavarapu2014nitsche} have pointed out that the penalty method does not converge for this problem unless the penalty parameter is chosen very carefully through trial and error, and that even if it converges, the resulting solution exhibits severe oscillations in normal and tangential tractions. 
To examine whether the use of the smoothed friction law can remedy this problem, we repeat the same problem using the smoothed friction law and the penalty method.
To isolate the effect of normal contact behavior, we use the two different methods only for the tangential traction, while treating the normal traction commonly by the barrier approach.
We assign the tangential penalty parameter as $\alpha_{\ct}=1000E_{\text{hard}}$, without an attempt to find a parameter that gives a converging (but perhaps erroneous) solution.

\begin{figure}[h!]
    \centering
    \subfloat[$E_{\mathrm{hard}}/E_{\mathrm{soft}}=10^1$]{\includegraphics[width=0.45\textwidth]{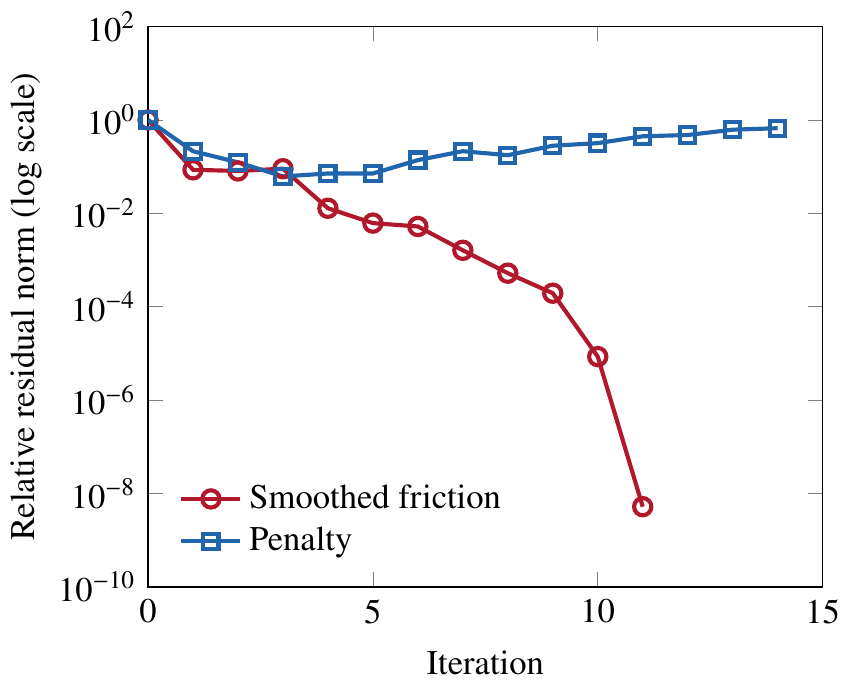}}
    \subfloat[$E_{\mathrm{hard}}/E_{\mathrm{soft}}=10^7$]{\includegraphics[width=0.45\textwidth]{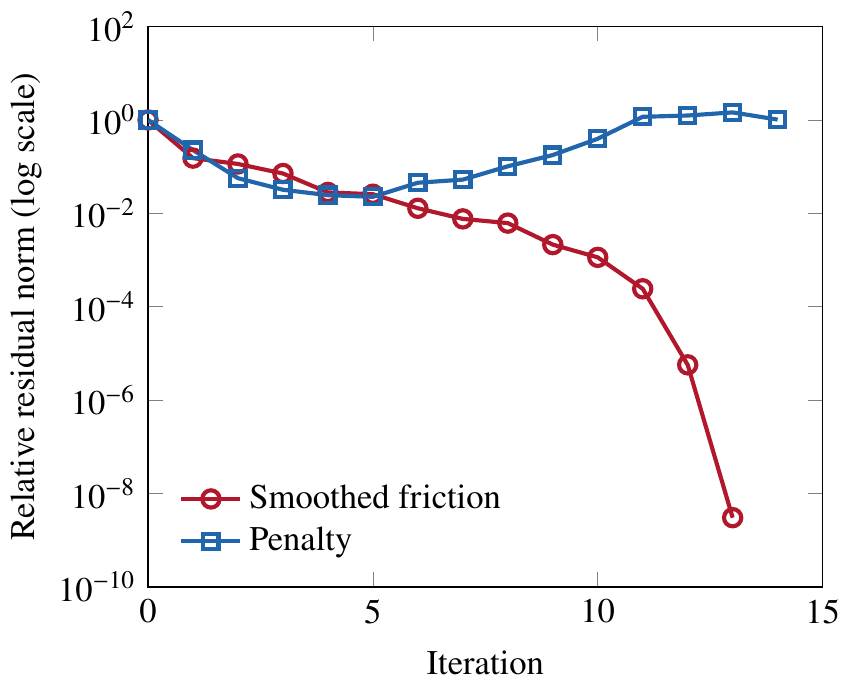}}$\quad$
    \caption{Sliding between two blocks with stiffness contrast: Newton convergence profiles of the smoothed friction law and the penalty method.}
    \label{fig:stiffness-contrast-convergence}
\end{figure}

Figure~\ref{fig:stiffness-contrast-convergence} shows the Newton convergence profiles of the smoothed friction law and the penalty method.
Being consistent with what reported by Annavarapu \etal~\cite{annavarapu2014nitsche}, the penalty method does not converge in both cases, despite the use of a direct linear solver.
While not presented for brevity, similar non-convergence occurs when the normal contact is treated by the penalty method instead of the barrier approach.
Conversely, when the smoothed friction law is used, the simulation does converge in both cases, showing optimal (asymptotically quadratic) rates of convergence.
From the numerical viewpoint, the main difference between the two methods is that the tangent operator of the tangential traction vector, $\pd\bm{t}_{\ct}/\pd u_{\ct}$, changes discontinuously during a stick--slip transition in the penalty method, whereas it varies continuously with the slip displacement in the smoothed friction law.
It can be seen that this difference provides significant numerical robustness to problems where stick and slip modes are mixed.

\begin{figure}[h!]
    \centering
    \includegraphics[width=1.0\textwidth]{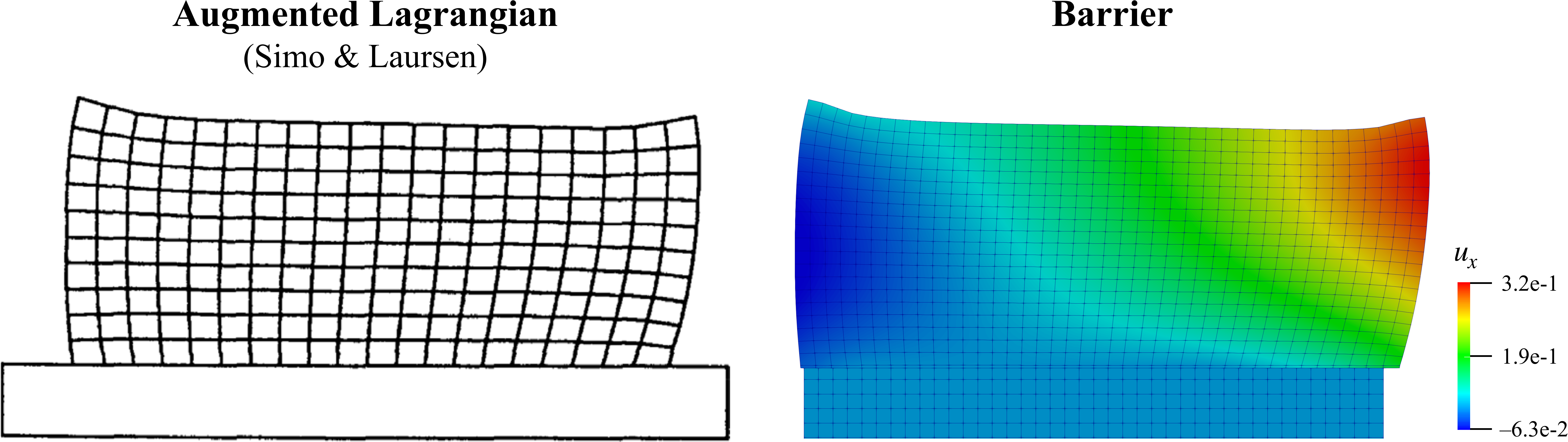}
    \caption{Sliding between two blocks with stiffness contrast: comparison of deformed geometries obtained by Simo and Laursen~\cite{simo1992augmented} and the current method. The barrier result also shows the $x$-displacement field in meter.}
    \label{fig:stiffness-contrast-deformed-geometry}
\end{figure}

\begin{figure}[h!]
    \centering
    \subfloat[Contant pressure]{\includegraphics[width=0.45\textwidth]{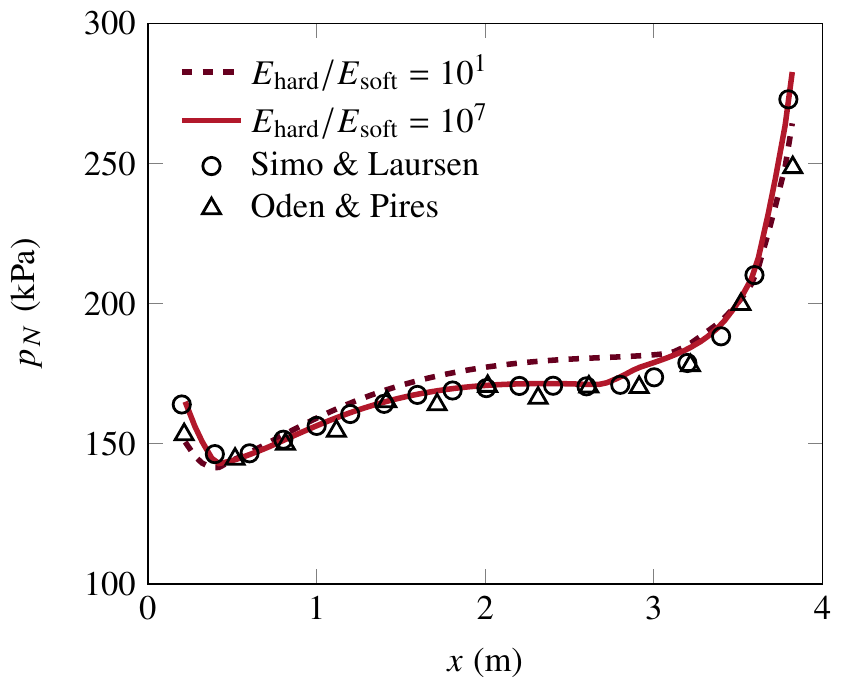}}$\quad$
    \subfloat[Tangential stress]{\includegraphics[width=0.45\textwidth]{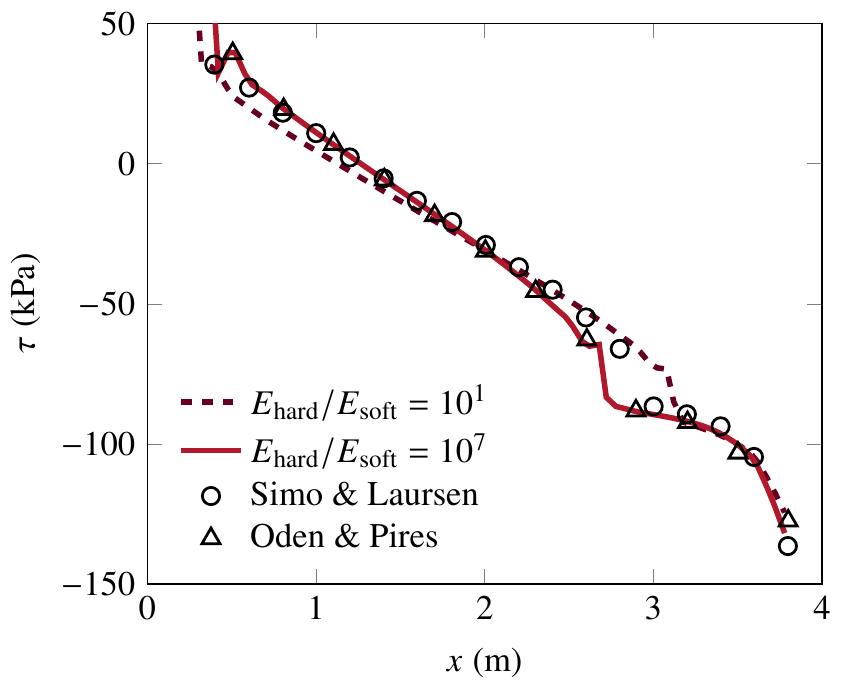}}
    \caption{Sliding between two blocks with stiffness contrast: contact pressures and tangential stresses. The reference results from Oden and Pires~\cite{oden1984algorithms} and Simo and Laursen~\cite{simo1992augmented} are obtained when the lower block is rigid.}
    \label{fig:stiffness-contrast-pN-tau}
\end{figure}

We then verify the numerical solutions obtained by the smoothed friction law in Figs.~\ref{fig:stiffness-contrast-deformed-geometry} and~\ref{fig:stiffness-contrast-pN-tau}.
As shown in Fig.~\ref{fig:stiffness-contrast-deformed-geometry}, the deformed geometry of our numerical solution to the high stiffness contrast case is virtually identical to that obtained by Simo and Laursen~\cite{simo1992augmented} using an augmented Lagrangian method along with the standard (non-embedded) finite element method.
Figure~\ref{fig:stiffness-contrast-pN-tau} also demonstrates that the contact pressures and tangential stresses along the interface do not show any oscillations, unlike the penalty solutions (and the weighted Nitsche solutions) in Annavarapu \etal~\cite{annavarapu2014nitsche}.
Note, however, that the proposed method is still subjected to severe pressure oscillations.
In the next example, we demonstrate a case where the barrier method shows severe oscillations and how to suppress the oscillations using the averaged surface integration scheme.

\subsection{Separation of a circular inclusion}
\begin{figure}[h!]
    \centering
    \includegraphics[width=0.5\textwidth]{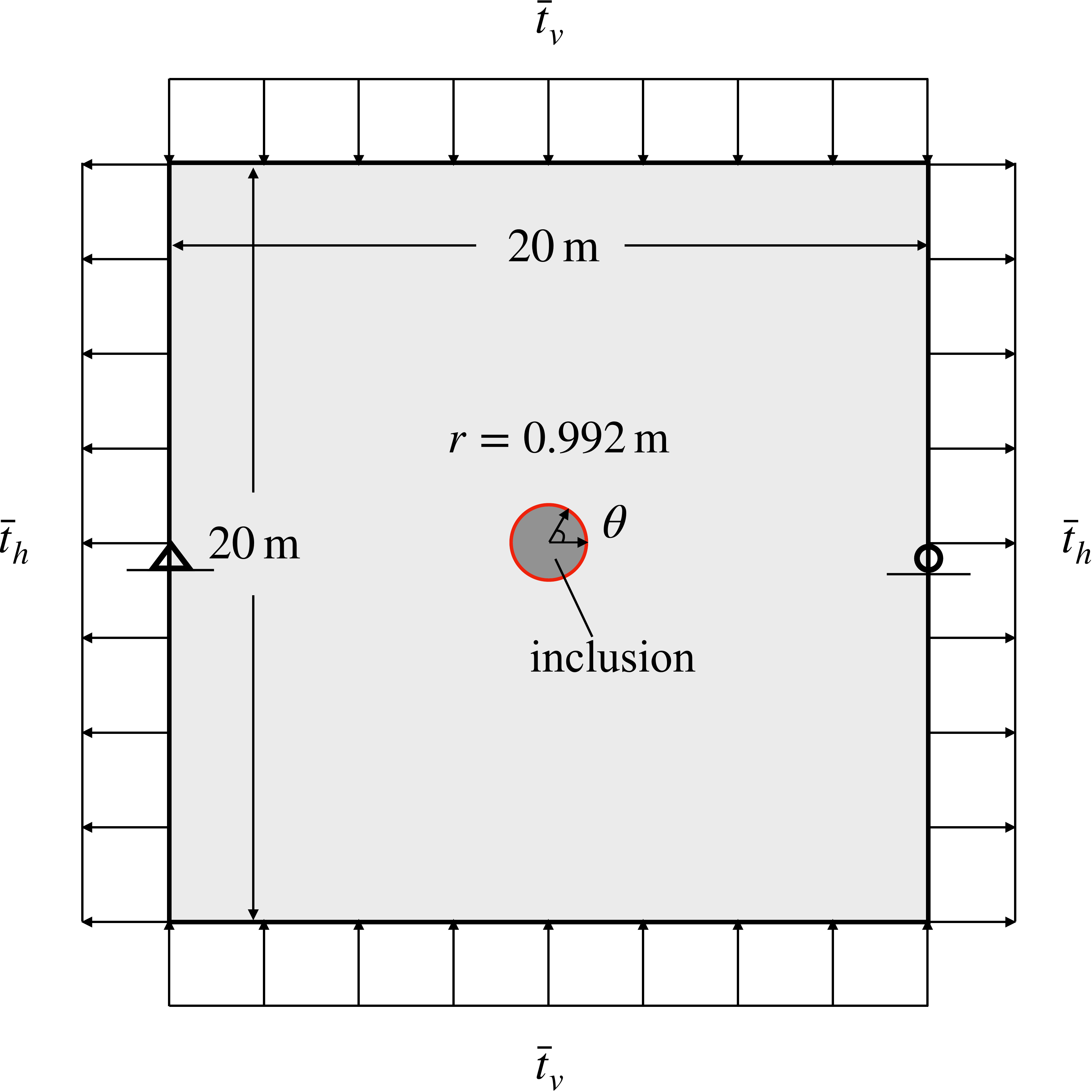}
    \caption{Separation of a circular inclusion: problem geometry and boundary conditions. $\bar{t}$ denotes the magnitude of the external traction.}
    \label{fig:inclusion-setup}
\end{figure}

\revised{The purpose of our final example is twofold: (i) to investigate the performance of the averaged surface integration scheme to alleviate traction oscillations, and (ii) to demonstrate the robustness of the barrier treatment when the contact condition changes with load.}
To this end, we consider the separation of a circular inclusion in an elastic domain due to to compressive and tensile tractions, which is depicted in Fig.~\ref{fig:inclusion-setup}.
Keer \etal~\cite{keer1973separation} derived an analytical solution to this problem in the case of that the surrounding domain is infinitely large, the inclusion and the surrounding material have the same stiffness, and the interface between the inclusion and the surrounding material is frictionless.
So we first consider this particular case to verify our solution against the analytical solution. 
According to the analytical solution, some part of the interface (roughly, $\theta < 55^{\circ}$) would be separated, while the rest is in contact and shows positive contact pressures.

To simulate the case that permits an analytical solution, we assign the Young's modulus and Poisson's ratio of the inclusion and the surrounding material as $E=1000$ kPa and $\nu=0$, and the boundary tractions $\bar{t}_{h}=\bar{t}_{v}=10$ kPa.
To investigate the convergence of the averaged integration scheme, we use three different meshes comprised of mono-sized square elements: Mesh 1 ($h=0.2$ m), Mesh 2 ($h=0.1$ m), and Mesh 3 ($h=0.05$ m). 
In each mesh, the center and its adjacent nodes inside the inclusion are supported by horizontal rollers to prevent the rotation of the inclusion.
We compute the barrier stiffness parameter $\kappa$ approximating the value of $(p_{N,0})_{\text{optimal}}$ to be the external traction.

\begin{figure}[h!]
    \centering
    \subfloat[Standard integration scheme]{\includegraphics[width=0.45\textwidth]{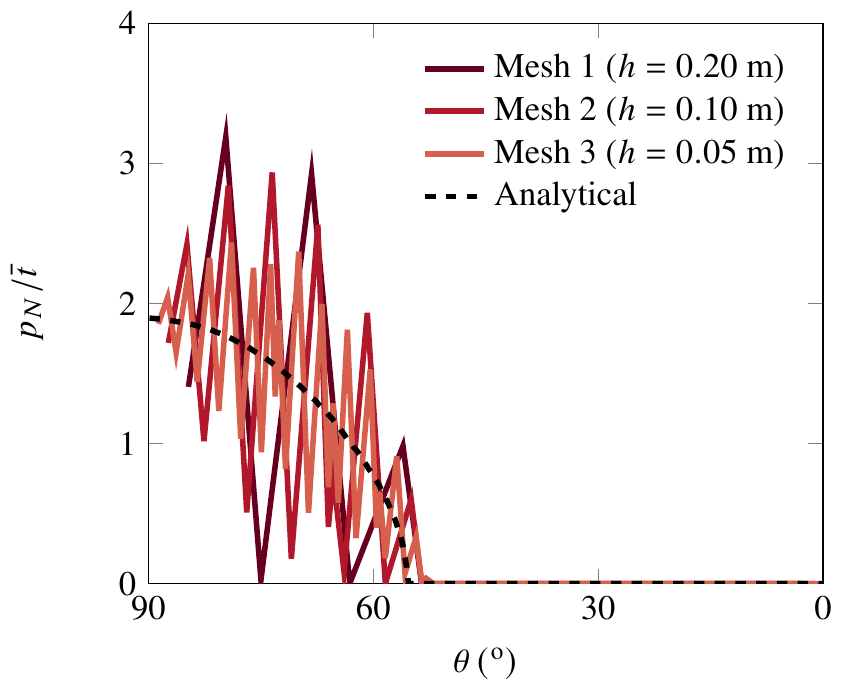}}$\quad$
    \subfloat[Averaged integration scheme]{\includegraphics[width=0.45\textwidth]{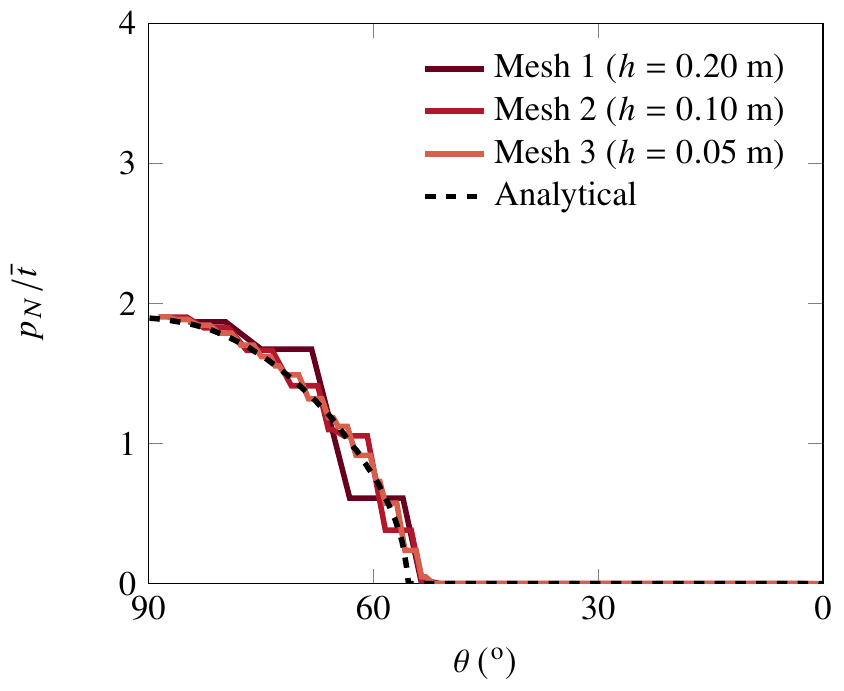}}
    \caption{Separation of a circular inclusion: contact pressure solutions obtained by (a) the standard surface integration scheme and (b) the averaged surface integration scheme.}
    \label{fig:circular-integration-default}
\end{figure}

\begin{figure}[h!]
    \centering
    \includegraphics[width=0.45\textwidth]{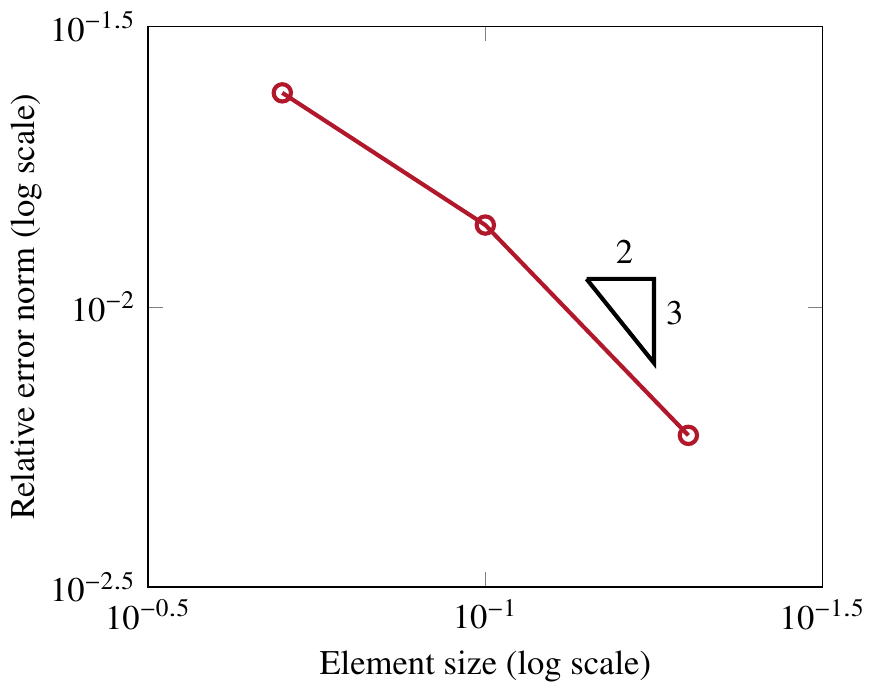}
    \caption{\revised{Separation of a circular inclusion: convergence of the contact pressure solutions with mesh refinement.}}
    \label{fig:circular-integration-convergence}
\end{figure}

Figure~\ref{fig:circular-integration-default} shows contact pressure solutions along the interface obtained by the standard integration scheme and the averaged integration scheme, comparing them with the analytical solution (the pressures are normalized by $\bar{t}\equiv\bar{t}_{v}=\bar{t}_{h}$).
It can be seen that the solutions obtained by the standard scheme show severe oscillations which are not well suppressed by mesh refinement.
Conversely, those obtained by the averaged scheme do not show any oscillations and converge well to the analytical solution as the element size decreases.
Notably, the results are slightly less oscillatory than the numerical solutions obtained by the augmented Lagrangian method in Belytschko \etal~\cite{belytschko2001arbitrary}.
\revised{In Fig.~\ref{fig:circular-integration-convergence}, we also plot the relative error norms of the pressure field solutions with the element size.
It can be seen that the numerical solutions show a good asymptotic convergence rate, which is comparable to that of Lagrange multiplier solutions in Kim \etal~\cite{kim2007mortared}.}
These results demonstrate that the averaged integration scheme is an effective and consistent method to suppress traction oscillations.
It is also noted that while not presented for brevity, we observe the same from the penalty solutions to this problem.

Having verified the averaged integration scheme against the analytical solution, we extend the problem to more challenging cases, by introducing (i) stiffness contrast between the inclusion and the surrounding material, and then (ii) friction along the interface.
For the former extension, we increase the Young's modulus of the inclusion by 100 times. For the latter, we set the friction coefficient of the interface as $\mu = 0.3$.
The resulting problem may be akin to a typical soil--pipe interaction problem.
We use Mesh 3 (the finest mesh) for both cases.

\begin{figure}[h!]
    \centering
    \subfloat[With stiffness contrast]{\includegraphics[width=0.45\textwidth]{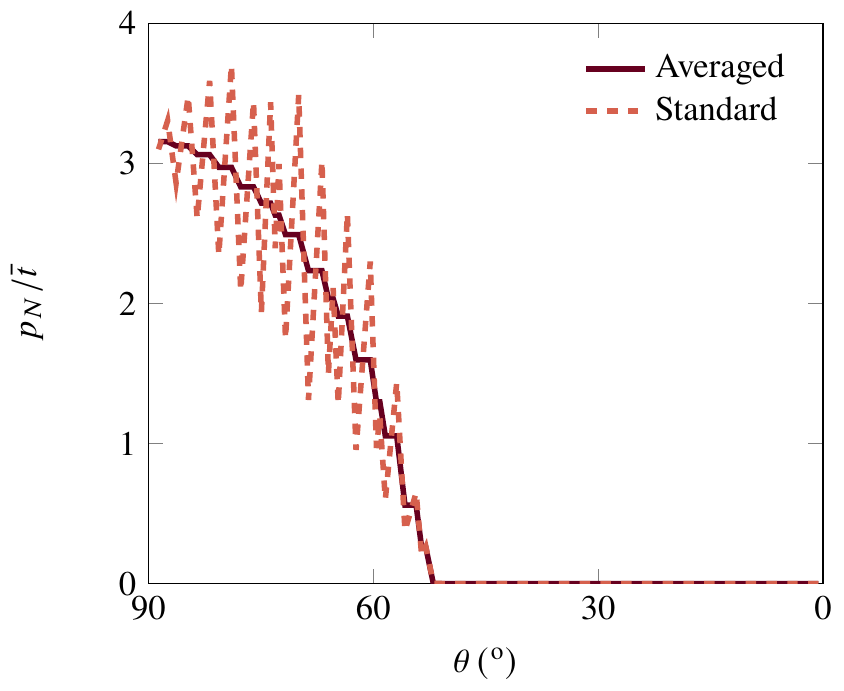}}$\quad$
    \subfloat[With stiffness contrast and friction]{\includegraphics[width=0.45\textwidth]{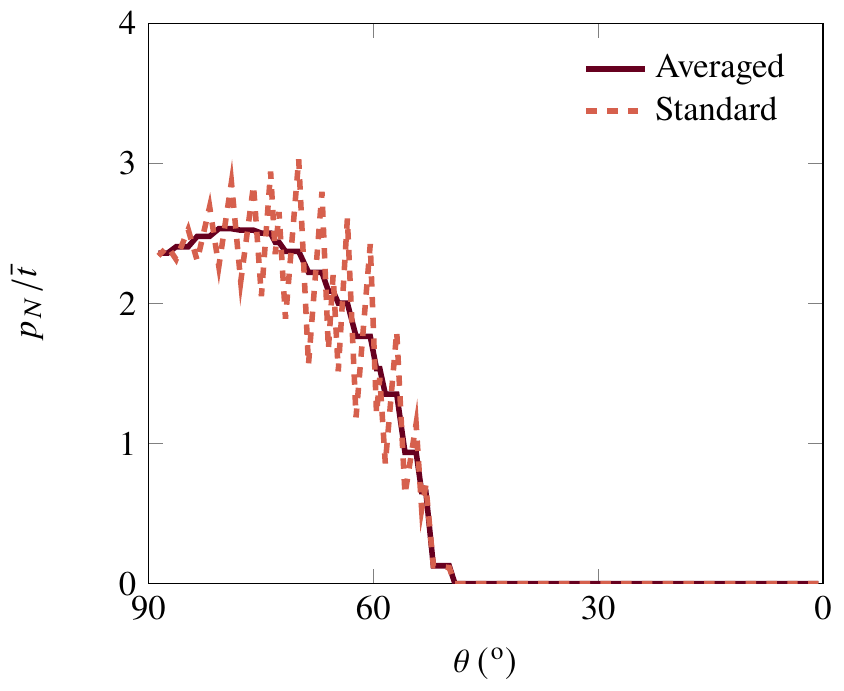}}
    \caption{Separation of a circular inclusion (a) with stiffness contrast and (b) with stiffness contrast and friction: contact pressure solutions obtained by the standard and averaged surface integration schemes.}
    \label{fig:circular-integration-extended}
\end{figure}

Figure~\ref{fig:circular-integration-extended} presents the contact pressure solutions in the two extended cases.
We confirm that the averaged integration scheme continues to perform well for these more challenging cases.
So it can be concluded that the combination of the barrier method and the averaged surface integration scheme is robust for highly challenging frictional contact problems, while it can be utilized at appreciably low computational cost.  

\begin{figure}[h!]
    \centering
    \includegraphics[width=1.0\textwidth]{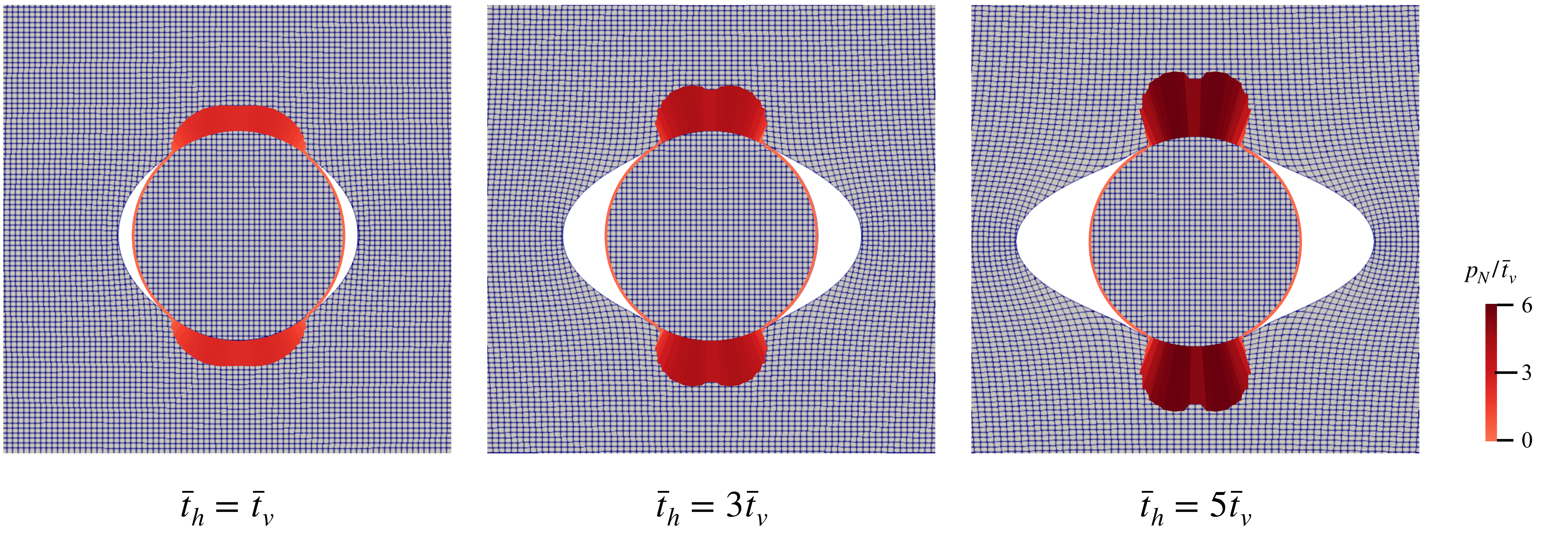}
    \caption{\revised{Separation of a circular inclusion: deformed meshes and normalized contact pressures ($p_{\cn}/\bar{t}$) during an increase in the horizontal traction. Displacements are exaggerated by a factor of 5.}}
    \label{fig:circular-integration-deformed-mesh}
\end{figure}

\revised{Lastly, to investigate the robustness of the barrier treatment under evolving contact conditions, we further extend the problem by increasing the magnitude of the tensile traction in the horizontal direction, $\bar{t}_{h}$, from $\bar{t}_{h}=\bar{t}_{v}$ to $\bar{t}_{h}=5\bar{t}_{v}$.
Figure~\ref{fig:circular-integration-deformed-mesh} presents the deformed meshes during this course of loading, superimposing the distribution of $p_{\cn}/\bar{t}_{v}$ along the interface calculated by the averaged integration scheme.
It can be seen that as the domain is pulled more in the horizontal direction, the contacting area becomes narrower and hence the contact pressure increases.
The results affirm that the barrier method continues to provide robust numerical solutions under evolving contact conditions.}

% SECTION
% ==============================================================================
\section{Closure}
We have developed the first barrier method for frictional contact on interfaces embedded in finite elements.
The method derives the contact pressure from a smooth barrier energy function that satisfies the contact constraints.
Also applied in the proposed method is a smoothed friction law where the stick--slip transition is approximated by a continuous function of the slip displacement.
The resulting formulation has the same mathematical structure as the penalty method, and hence it can be well cast into an existing penalty-based XFEM.
Drawing on this similarity, we have also devised the averaged surface integration scheme which can alleviate traction oscillations in the barrier method and other penalty-type methods.
Last but not least, we have derived values for the parameters of the barrier method specialized to embedded interfaces, so that the method can be practically used without any parameter tuning.

Through various numerical examples, we have demonstrated that the proposed barrier method features a highly competitive combination of robustness, accuracy, and efficiency. 
The barrier method has practically the same computational cost as the penalty method, which has been popular in practice due to its unparalleled efficiency.
Unlike the penalty method, however, the barrier method guarantees satisfaction of the non-penetration constraint with controllable accuracy. This feature is particularly desirable when inter-penetration should be strictly avoided, like coupled multiphysical problems in which the fracture opening should be non-negative.
The proposed method is also highly robust to common numerical issues such as non-convergence and erroneous traction oscillations.
Particularly, the barrier method aided by the averaged integration scheme has provided convergent and non-oscillatory solutions to all problems we have tested (besides those presented in this paper).
For these reasons, we believe that the proposed method is one of the most attractive options for handling frictional sliding in embedded finite element methods.

\revised{
It is believed that the proposed barrier method will continue to perform well for frictional contact under more complex conditions such as three-dimensional geometry, evolving discontinuity, and sophisticated friction laws. 
The main reason is that the proposed method essentially provides an alternative expression to the pressure--displacement relationship of the classic penalty method, which has been well applied to embedded interfaces in more challenging settings (\eg~\cite{liu2009extended,liu2013extended,prevost2016faults}).
The method may also be extended to handle frictional contact on interfaces embedded in other types of numerical methods (\eg~\cite{liang2022shear}).
Note, however, that extending the current formulation to accommodate these additional complexities entail lots of effort on other aspects (\eg~algorithms for handling three-dimensional and/or evolving interface geometry in embedded elements). We thus leave such an extension as a topic of future research.}

% ACKNOWLEDGEMENT
% ------------------------------------------------------------------------------
\section*{Acknowledgments}
The authors wish to thank the two reviewers for their expert comments.
Portions of the work were supported by the Research Grants Council of Hong Kong through Projects 17201419 and 27205918. 
Financial support from KAIST is also acknowledged.

% APPENDICES
% ==============================================================================
\appendix

\section{Optimal value of initial normal stiffness}
\label{appendix:A}

This appendix describes how to find the (approximately) optimal value of initial normal stiffness, $k_{N,0}$, that gives the best possible condition of the Jacobian matrix~\eqref{eq:jacobian}.
To begin, we recall that $k_{N, 0}$ enters the (2,2) block of the Jacobian matrix.
The element-wise contribution of this block can be written as given by
\begin{linenomath}
\begin{align}
	\left[ \frac{\pd \bm{\mathcal{R}}^{a}}{\pd \bm{a}} \right]^{e}
	&= \underbrace{\int_{\Omega \setminus \Gamma} \bm{B}_{I}^{\sf T}\bm{D}\bm{B}_J \,\dd V^{e}}_{=: \jump{\bm{K}}_{IJ}^e}
	\,+\, \underbrace{ \int_{\Gamma} N_I[k_{\cn,0} (\bm{n} \otimes \bm{n})]N_J \, \dd A}_{=: \bm{E}_{IJ}^e}.
    \label{eq:jacobian-2-2-block}
\end{align}
\end{linenomath}
where $\bm{B} := \symgrad N$ is the strain--displacement matrix, 
and $\bm{D}$ is the stiffness matrix of the bulk material.

For simplicity, we assume that the matrix~\eqref{eq:jacobian-2-2-block} is best conditioned when the max norms of $\jump{\bm{K}}_{IJ}^e$ and $\bm{E}_{IJ}^e$ are equal.
To calculate the max norm of $\jump{\bm{K}}_{IJ}^e$, we recall that the max norm of the $\bm{B}$ matrix is $1/h$, where $h$ is the element size.
The max norm of $\bm{D}$ may be approximated to be $E$.
Lastly, the integration volume is $h^{\rm dim}$, where dim denotes the spatial dimension.
Therefore, the max norm of $\jump{\bm{K}}_{IJ}^e$ can be approximated as
\begin{equation}
	\|\jump{\bm{K}}_{IJ}^e\|_{\rm max} \approx \sum_Q \left(\dfrac{1}{h} E \dfrac{1}{h}\right) h^{\rm dim} = \sum_Q Eh^{{\rm dim}-2}.
	\label{eq:mag_K}
\end{equation}
Similarly, the max norm of $\bm{E}_{IJ}^e$ can be approximated as (note that $k_{N, 0}>0$ by definition)
\begin{equation}
	\|\bm{E}_{IJ}^e\|_{\rm max} \approx \sum_Q k_{N,0}\, h^{{\rm dim}-1}. 
	\label{eq:mag_E} 
\end{equation} 
The optimal value of the initial barrier stiffness can be obtained by equating Eqs.~\eqref{eq:mag_K} and~\eqref{eq:mag_E}, which gives
\begin{equation}
	(k_{N, 0})_{\mathrm{optimal}} = E/h.
\end{equation}

% REFERENCES
% ------------------------------------------------------------------------------
% \section*{References}
\bibliography{references}

\begin{thebibliography}{10}
\expandafter\ifx\csname url\endcsname\relax
  \def\url#1{\texttt{#1}}\fi
\expandafter\ifx\csname urlprefix\endcsname\relax\def\urlprefix{URL }\fi
\expandafter\ifx\csname href\endcsname\relax
  \def\href#1#2{#2} \def\path#1{#1}\fi

\bibitem{simo1990class}
J.~C. Simo, M.~Rifai, A class of mixed assumed strain methods and the method of
  incompatible modes, International Journal for Numerical Methods in
  Engineering 29~(8) (1990) 1595--1638.

\bibitem{moes1999finite}
N.~Mo{\"e}s, J.~E. Dolbow, T.~Belytschko, A finite element method for crack
  growth without remeshing, International Journal for Numerical Methods in
  Engineering 46~(1) (1999) 131--150.

\bibitem{dolbow2001extended}
J.~E. Dolbow, N.~Mo{\"e}s, T.~Belytschko, An extended finite element method for
  modeling crack growth with frictional contact, Computer Methods in Applied
  Mechanics and Engineering 190~(51-52) (2001) 6825--6846.

\bibitem{khoei2007enriched}
A.~Khoei, M.~Nikbakht, An enriched finite element algorithm for numerical
  computation of contact friction problems, International Journal of Mechanical
  Sciences 49~(2) (2007) 183--199.

\bibitem{liu2008contact}
F.~Liu, R.~I. Borja, A contact algorithm for frictional crack propagation with
  the extended finite element method, International Journal for Numerical
  Methods in Engineering 76~(10) (2008) 1489--1512.

\bibitem{nistor2009xfem}
I.~Nistor, M.~L.~E. Guiton, P.~Massin, N.~Mo{\"e}s, S.~Geniaut, {An X-FEM
  approach for large sliding contact along discontinuities}, International
  Journal for Numerical Methods in Engineering 78~(12) (2009) 1407--1435.

\bibitem{elguedj2007mixed}
T.~Elguedj, A.~Gravouil, A.~Combescure, A mixed augmented lagrangian-extended
  finite element method for modelling elastic--plastic fatigue crack growth
  with unilateral contact, International Journal for Numerical Methods in
  Engineering 71~(13) (2007) 1569--1597.

\bibitem{bechet2009stable}
{\'E}.~B{\'e}chet, N.~Mo{\"e}s, B.~Wohlmuth, {A stable Lagrange multiplier
  space for stiff interface conditions within the extended finite element
  method}, International Journal for Numerical Methods in Engineering 78~(8)
  (2009) 931--954.

\bibitem{liu2010stabilized}
F.~Liu, R.~I. Borja, Stabilized low-order finite elements for frictional
  contact with the extended finite element method, Computer Methods in Applied
  Mechanics and Engineering 199~(37-40) (2010) 2456--2471.

\bibitem{laursen2003computational}
T.~A. Laursen, Computational Contact and Impact Mechanics, Springer, 2003.

\bibitem{wriggers2006computational}
P.~Wriggers, Computational Contact Mechanics, Springer, 2006.

\bibitem{kim2007mortared}
T.~Y. Kim, J.~E. Dolbow, T.~A. Laursen, A mortared finite element method for
  frictional contact on arbitrary interfaces, Computational Mechanics 39~(3)
  (2007) 223--235.

\bibitem{ji2004strategies}
H.~Ji, J.~E. Dolbow, On strategies for enforcing interfacial constraints and
  evaluating jump conditions with the extended finite element method,
  International Journal for Numerical Methods in Engineering 61~(14) (2004)
  2508--2535.

\bibitem{alart1991mixed}
P.~Alart, A.~Curnier, A mixed formulation for frictional contact problems prone
  to newton like solution methods, Computer methods in applied mechanics and
  engineering 92~(3) (1991) 353--375.

\bibitem{simo1992augmented}
J.~C. Simo, T.~A. Laursen, {An augmented Lagrangian treatment of contact
  problems involving friction}, Computers \& Structures 42~(1) (1992) 97--116.

\bibitem{annavarapu2013nitsche}
C.~Annavarapu, M.~Hautefeuille, J.~E. Dolbow, {A Nitsche stabilized finite
  element method for frictional sliding on embedded interfaces. Part II:
  Intersecting interfaces}, Computer Methods in Applied Mechanics and
  Engineering 267 (2013) 318--341.

\bibitem{annavarapu2014nitsche}
C.~Annavarapu, M.~Hautefeuille, J.~E. Dolbow, {A Nitsche stabilized finite
  element method for frictional sliding on embedded interfaces. Part I: single
  interface}, Computer Methods in Applied Mechanics and Engineering 268 (2014)
  417--436.

\bibitem{annavarapu2015weighted}
C.~Annavarapu, R.~R. Settgast, S.~M. Johnson, P.~Fu, E.~B. Herbold, {A weighted
  Nitsche stabilized method for small-sliding contact on frictional surfaces},
  Computer Methods in Applied Mechanics and Engineering 283 (2015) 763--781.

\bibitem{fei2020phase-a}
F.~Fei, J.~Choo, A phase-field method for modeling cracks with frictional
  contact, International Journal for Numerical Methods in Engineering 121~(4)
  (2020) 740--762.

\bibitem{fei2020phase-b}
F.~Fei, J.~Choo, A phase-field model of frictional shear fracture in geologic
  materials, Computer Methods in Applied Mechanics and Engineering 369 (2020)
  113265.

\bibitem{fei2021double}
F.~Fei, J.~Choo, Double-phase-field formulation for mixed-mode fracture in
  rocks, Computer Methods in Applied Mechanics and Engineering 376 (2021)
  113655.

\bibitem{fei2021phase}
F.~Fei, J.~Choo, C.~Liu, J.~A. White, Phase-field modeling of rock fractures
  with roughness, arXiv preprint arXiv:2105.14663 (2021).

\bibitem{kane1999finite}
C.~Kane, E.~A. Repetto, M.~Ortiz, J.~E. Marsden, Finite element analysis of
  nonsmooth contact, Computer Methods in Applied Mechanics and Engineering
  180~(1-2) (1999) 1--26.

\bibitem{pandolfi2002time}
A.~Pandolfi, C.~Kane, J.~E. Marsden, M.~Ortiz, Time-discretized variational
  formulation of non-smooth frictional contact, International Journal for
  Numerical Methods in Engineering 53~(8) (2002) 1801--1829.

\bibitem{li2020incremental}
M.~Li, Z.~Ferguson, T.~Schneider, T.~Langlois, D.~Zorin, D.~Panozzo, C.~Jiang,
  D.~M. Kaufman, Incremental potential contact: Intersection-and
  inversion-free, large-deformation dynamics, ACM Transactions on Graphics
  39~(4) (2020).

\bibitem{Li2021CIPC}
M.~Li, D.~M. Kaufman, C.~Jiang, Codimensional incremental potential contact,
  ACM Transactions on Graphics 40~(4) (2021).

\bibitem{Ferguson:2021:RigidIPC}
Z.~Ferguson, M.~Li, T.~Schneider, F.~Gil-Ureta, T.~Langlois, C.~Jiang,
  D.~Zorin, D.~M. Kaufman, D.~Panozzo, Intersection-free rigid body dynamics,
  ACM Transactions on Graphics 40~(4) (2021).

\bibitem{Lan2021MIPC}
L.~Lan, Y.~Yang, D.~M. Kaufman, J.~Yao, M.~Li, C.~Jiang, Medial {IPC}:
  Accelerated incremental potential contact with medial elastics, ACM
  Transactions on Graphics 40~(4) (2021).

\bibitem{liu2017stabilized}
F.~Liu, P.~Gordon, H.~Meier, D.~Valiveti, A stabilized extended finite element
  framework for hydraulic fracturing simulations, International Journal for
  Numerical and Analytical Methods in Geomechanics 41~(5) (2017) 654--681.

\bibitem{choo2018cracking}
J.~Choo, W.~Sun, Cracking and damage from crystallization in pores: Coupled
  chemo-hydro-mechanics and phase-field modeling, Computer Methods in Applied
  Mechanics and Engineering 335 (2018) 347--379.

\bibitem{liu2020modeling}
F.~Liu, Modeling hydraulic fracture propagation in permeable media with an
  embedded strong discontinuity approach, International Journal for Numerical
  and Analytical Methods in Geomechanics 44~(12) (2020) 1634--1655.

\bibitem{boyd2004convex}
S.~Boyd, L.~Vandenberghe, Convex Optimization, Cambridge University Press,
  2004.

\bibitem{wriggers1990finite}
P.~Wriggers, T.~V. Van, E.~Stein, Finite element formulation of large
  deformation impact-contact problems with friction, Computers \& Structures
  37~(3) (1990) 319--331.

\bibitem{sandeep2019experimental}
C.~Sandeep, K.~Senetakis, An experimental investigation of the microslip
  displacement of geological materials, Computers and Geotechnics 107 (2019)
  55--67.

\bibitem{cusini2021simulation}
M.~Cusini, J.~A. White, N.~Castelletto, R.~R. Settgast, Simulation of coupled
  multiphase flow and geomechanics in porous media with embedded discrete
  fractures, International Journal for Numerical and Analytical Methods in
  Geomechanics 45~(5) (2021) 563--584.

\bibitem{borja2008assumed}
R.~I. Borja, Assumed enhanced strain and the extended finite element methods: A
  unification of concepts, Computer Methods in Applied Mechanics and
  Engineering 197~(33-40) (2008) 2789--2803.

\bibitem{bochev2006stabilization}
P.~B. Bochev, C.~R. Dohrmann, M.~D. Gunzburger, Stabilization of low-order
  mixed finite elements for the stokes equations, SIAM Journal on Numerical
  Analysis 44~(1) (2006) 82--101.

\bibitem{choo2015stabilized}
J.~Choo, R.~I. Borja, Stabilized mixed finite elements for deformable porous
  media with double porosity, Computer Methods in Applied Mechanics and
  Engineering 293 (2015) 131--154.

\bibitem{choo2019stabilized}
J.~Choo, {Stabilized mixed continuous/enriched Galerkin formulations for
  locally mass conservative poromechanics}, Computer Methods in Applied
  Mechanics and Engineering 357 (2019) 112568.

\bibitem{zhao2020stabilized}
Y.~Zhao, J.~Choo, Stabilized material point methods for coupled large
  deformation and fluid flow in porous materials, Computer Methods in Applied
  Mechanics and Engineering 362 (2020) 112742.

\bibitem{svenning2016weak}
E.~Svenning, A weak penalty formulation remedying traction oscillations in
  interface elements, Computer Methods in Applied Mechanics and Engineering 310
  (2016) 460--474.

\bibitem{liu2010finite}
F.~Liu, R.~I. Borja, Finite deformation formulation for embedded frictional
  crack with the extended finite element method, International Journal for
  Numerical Methods in Engineering 82~(6) (2010) 773--804.

\bibitem{arndt2021deal}
D.~Arndt, W.~Bangerth, D.~Davydov, T.~Heister, L.~Heltai, M.~Kronbichler,
  M.~Maier, J.-P. Pelteret, B.~Turcksin, D.~Wells, {The deal. II finite element
  library: Design, features, and insights}, Computers \& Mathematics with
  Applications 81 (2021) 407--422.

\bibitem{choo2018enriched}
J.~Choo, S.~Lee, Enriched {G}alerkin finite elements for coupled poromechanics
  with local mass conservation, Computer Methods in Applied Mechanics and
  Engineering 341 (2018) 311--332.

\bibitem{choo2018large}
J.~Choo, Large deformation poromechanics with local mass conservation: An
  enriched {G}alerkin finite element framework, International Journal for
  Numerical Methods in Engineering 116 (2018) 66--90.

\bibitem{oden1984algorithms}
J.~T. Oden, E.~B. Pires, Algorithms and numerical results for finite element
  approximations of contact problems with non-classical friction laws,
  Computers \& Structures 19~(1-2) (1984) 137--147.

\bibitem{keer1973separation}
L.~M. Keer, J.~Dundurs, K.~Kiattikomol, Separation of a smooth circular
  inclusion from a matrix, International Journal of Engineering Science 11~(11)
  (1973) 1221--1233.

\bibitem{belytschko2001arbitrary}
T.~Belytschko, N.~Mo{\"e}s, S.~Usui, C.~Parimi, Arbitrary discontinuities in
  finite elements, International Journal for Numerical Methods in Engineering
  50~(4) (2001) 993--1013.

\bibitem{liu2009extended}
F.~Liu, R.~I. Borja, An extended finite element framework for slow-rate
  frictional faulting with bulk plasticity and variable friction, International
  Journal for Numerical and Analytical Methods in Geomechanics 33~(13) (2009)
  1535--1560.

\bibitem{liu2013extended}
F.~Liu, R.~I. Borja, Extended finite element framework for fault rupture
  dynamics including bulk plasticity, International Journal for Numerical and
  Analytical Methods in Geomechanics 37~(18) (2013) 3087--3111.

\bibitem{prevost2016faults}
J.~H. Pr\'{e}vost, N.~Sukumar, Faults simulations for three-dimensional
  reservoir-geomechanical models with the extended finite element method,
  Journal of the Mechanics and Physics of Solids 86 (2016) 1--18.

\bibitem{liang2022shear}
Y.~Liang, B.~Chandra, K.~Soga, {Shear band evolution and post-failure
  simulation by the extended material point method (XMPM) with localization
  detection and frictional self-contact}, Computer Methods in Applied Mechanics
  and Engineering 390 (2022) 114530.

\end{thebibliography}
% \printbibliography

\end{document}